\documentclass[11pt]{amsart}
\usepackage{float}
\usepackage[usenames]{color}
\usepackage{graphicx}
\usepackage{amssymb}
\usepackage{amscd}
\theoremstyle{plain}
\newtheorem{thm}{Theorem}[section]

\newtheorem{cor}[thm]{Corollary}

\theoremstyle{definition}

\newtheorem{rmk}[thm]{Remark}
\newtheorem{example}[thm]{Example}

\def\det{\mathop{\hbox {det}}\nolimits}

\def\im{\mathop{\hbox {Im}}\nolimits}

\def\ker{\mathop{\hbox{Ker}}\nolimits}

\newcommand{\fra}{\mathfrak{a}}
\newcommand{\frb}{\mathfrak{b}}

\newcommand{\frg}{\mathfrak{g}}
\newcommand{\frh}{\mathfrak{h}}
\newcommand{\fri}{\mathfrak{i}}
\newcommand{\frk}{\mathfrak{k}}

\newcommand{\frp}{\mathfrak{p}}

\newcommand{\frt}{\mathfrak{t}}

\newcommand{\bbC}{\mathbb{C}}

\newcommand{\bbN}{\mathbb{N}}

\newcommand{\be}{\begin {equation}}
\newcommand{\ee}{\end {equation}}
\newcommand{\bp}{\begin {proof}}
\newcommand{\ep}{\end {proof}}

\begin{document}

\title[Dirac series for complex $E_7$]
{Dirac series for complex $E_7$}

\author{Chao-Ping Dong}
\author{Kayue Daniel Wong}

\address[Dong]{School of Mathematical Sciences, Soochow University, Suzhou 215006,
P.~R.~China}
\email{chaopindong@163.com}

\address[Wong]{School of Science and Engineering, The Chinese University of Hong Kong, Shenzhen,
Guangdong 518172, P. R. China}
\email{kayue.wong@gmail.com}

\abstract{This paper classifies the Dirac series  for complex $E_7$. As applications, we verify a few conjectures raised in 2011, 2019 and 2020 for this exceptional Lie group. In particular, according to Conjecture 1.1 of Barbasch and Pand\v zi\'c \cite{BP}, our classification should be helpful for understanding the unitary dual of complex $E_7$.}
\endabstract

\subjclass[2010]{Primary 22E46.}

\keywords{Dirac cohomology, spin norm, unitary representation.}

\maketitle
\section{Introduction}
Although many results quoted in this paper hold in a much wider setting, we embark with a connected simple \emph{complex} Lie group $G$ with finite center since our main concern here is actually  complex $E_7$. Let $\theta$ be the Cartan involution. Then $K:=G^{\theta}$ is a maximal compact subgroup of $G$. Let $\frg_0$ (resp., $\frk_0$) be the Lie algebra of $G$ (resp., $K$).  Let $T$ be a maximal torus of $K$. Let $\fra_0=\sqrt{-1}\frt_0$. Put $A=\exp(\fra_0)$. Then $H:=TA$ is the fundamental Cartan subgroup of $G$. Let
$$
\frg_0=\frk_0 + \frp_0
$$
be the Cartan decomposition on the Lie algebra level. As usual, we drop the subscripts to denote the complexifications. Let $\langle \cdot, \cdot \rangle$ be the Killing form on $\frg$. Let $U(\frg)$ be the universal enveloping algebra of $\frg$ and let $C(\frp)$ be the Clifford algebra of $\frp$.

To give geometric construction of the discrete series, Parthasarathy introduced the Dirac operator $D\in U(\frg)\otimes C(\frp)$ whose square is a natural Laplacian on $G$ in 1972 \cite{P72}. See \eqref{D-square}. Let $\pi$ be an irreducible $(\frg, K)$ module. Choose a spin module  $S$ of $C(\frp)$. Let $\widetilde{K}$ be the subgroup of $K\times {\rm Spin \,\frp_0}$ consisting of the pairs $(k, s)$ such that ${\rm Ad}(k)=p(s)$, where ${\rm Ad}: K\to SO(\frp_0)$ is the adjoint action and $p: {\rm Spin}\, \frp_0 \to SO(\frp_0)$ is the spin double covering map.

The Dirac operator $D$ acts on $X\otimes S$.  Whenever $\pi$ is unitary, the operator $D$ is self-adjoint on $X\otimes S$; moreover, the eigenvalue of $D^2$ on any $\widetilde{K}$-type of  $X\otimes S$ is non-negative. This leads to Parthasarathy's Dirac operator inequality \cite{P80}. See \eqref{Dirac-inequality}. It is very effective in  detecting non-unitary modules and thus a valuable tool for studying the unitary dual $\widehat{G}$.

To sharpen the Dirac inequality, and to obtain a better understanding of $\widehat{G}$, Vogan introduced  Dirac cohomology \cite{Vog97} as the following $\widetilde{K}$-module:
\begin{equation}\label{def-Dirac-cohomology}
H_D(\pi):=\ker D/ \ker D \cap \im D.
\end{equation}
Vogan conjectured that  $H_D(\pi)$ (whenever nonzero) should reveal the infinitesimal character of $\pi$, suggesting that Dirac cohomology should be a finer invariant of $\pi$. This conjecture was proven by Huang and Pand\v zi\'c \cite{HP1} in 2002. See Theorem \ref{thm-HP}.  Since then, classifying $\widehat{G}^{\mathrm d}$---the set of all the equivalence classes of irreducible unitarizable $(\frg, K)$ modules with non-vanishing Dirac cohomology---became an interesting open problem.  As suggested by Huang, we call $\widehat{G}^{\mathrm d}$ the \emph{Dirac series} of $G$.

It is worth noting that when $\pi$ is unitary, $H_D(\pi)=\ker D=\ker D^2$. Thus $\pi$ is a Dirac series of $G$ if and only if zero is an eigenvalue of $D^2$ on $\pi \otimes S$. Therefore, Dirac series are exactly the members of $\widehat{G}$ on which the Dirac inequality is attained.

This paper aims to classify $\widehat{E}_7^{\mathrm{d}}$---the Dirac series of the connected and simply connected complex Lie group such that $\frg_0$ is of type $E_7$. Let us build up necessary notation to state our main result.  We identify
\begin{equation}\label{identifications}
\frg\cong \frg_{0} \oplus
\frg_0, \quad
\frh\cong \frh_{0} \oplus
\frh_0, \quad \frt\cong \{(x,-x) : x\in
\frh_{0} \}, \quad \fra \cong\{(x, x) : x\in
\frh_{0} \}.
\end{equation}
Denote the Weyl group of $\Delta(\frg_0, \frh_0)$ by $W$, which has identity element $e$ and longest element $w_0$.
Then $W(\Delta(\frg, \frh))$ is isomorphic to $W\times W$. By Theorem \ref{thm-Zh}, every irreducible admissible $(\frg, K)$ module is isomorphic to certain  $J(\lambda_L, \lambda_R)$, see Section \ref{sec-Zh}.

\begin{figure}[H]
\centering \scalebox{0.7}{\includegraphics{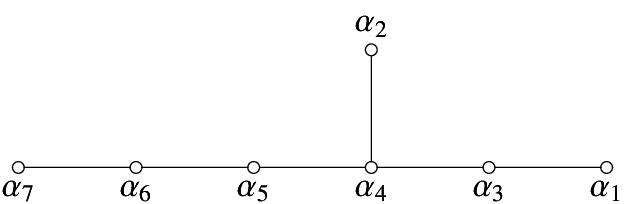}}
\caption{Dynkin diagram for $E_7$}
\label{E7-Dynkin}
\end{figure}

Let us fix the  simple roots of complex $E_7$ as in Fig.~\ref{E7-Dynkin}.
We enumerate the corresponding fundamental weights as $\varpi_1, \dots, \varpi_7$, and use them as a basis to express weights. That is, $\lambda=[\lambda_1, \cdots, \lambda_7]$ stands for the weight $\sum_{i=1}^{7}\lambda_i\varpi_i$. We say that $\lambda$ is \emph{integral} (resp., \emph{half-integral}) if each coordinate $\lambda_i$ (resp., $2\lambda_i$) is an integer. In this terminology, half sum of the positive roots in $\Delta^+(\frg_0, \frh_0)$ is expressed as $\rho(E_7)=[1, 1,1,1,1,1,1]$ and it is integral regular.
There are $10208$ involutions in the Weyl group $W(E_7)$ in total, among which $8479$ have the property that
\begin{equation}\label{def-Is}
I(s):=\{1\leq i\leq 7\mid s(\varpi_i)=\varpi_i\}
\end{equation}
is empty. As shown in Lemma 3.1 of \cite{DD}, the set $I(s)$ collects those $1\leq i\leq 7$ such that the simple reflection $s_i$ does \emph{not} occur in some (thus any) reduced expression of $s$. Let $L(s)$ be the Levi factor of the parabolic subgroup  of complex $E_7$ generated by the simple roots $\{\alpha_i\mid i \notin I(s)\}$.

Write $\{1, 2, 3, \cdots\}$ as $\mathbb{P}$ and $\{0, 1, 2,  \cdots\}$ as $\mathbb{N}$. For any fixed involution $s\in W(E_7)$, put
\begin{equation}\label{Lambda-s}
\Lambda(s):=\{[\lambda_1, \dots, \lambda_l] \mid 2\lambda_i\in\mathbb{P}, \lambda-s\lambda  \mbox { is a } \mathbb{N} \mbox{ combination of simple roots} \}.
\end{equation}
By page 5 of \cite{BP} and Section 2.6 of \cite{D2}, any representation in $\widehat{E}_7^{\mathrm{d}}$ must bear the form
\begin{equation}\label{can-repns}
J(\lambda, -s\lambda),
\end{equation}
where $s\in W(E_7)$ is an involution and $\lambda\in\Lambda(s)$.
This module has lowest $K$-type $\{\lambda+s\lambda\}$, the unique dominant weight to which $\lambda+s\lambda$ is conjugate under the action of $W(E_7)$. Whenever the set $I(s)$ is empty and that the module $J(\lambda, -s\lambda)$ has non-zero Dirac cohomology, we will call $J(\lambda, -s\lambda)$ a \emph{scattered representation}.  It is shown in Proposition 3.4 of \cite{DD} that there are finitely many scattered representations.  On the other hand, if the set $I(s)$ is non-empty, as shown in Proposition 3.5 of \cite{DD}, we can partition the modules $J(\lambda, -s\lambda)$ with non-zero Dirac cohomology (if they do exist) into finitely many parts, each bearing the form
\begin{equation}\label{string}
\{J(\lambda, -s\lambda)\mid  \lambda_i \mbox{ runs over } \frac{1}{2}\mathbb{P} \mbox{ for } i \in I(s) \mbox{ and } \lambda_j \mbox{ is fixed for } j \notin I(s) \}.
\end{equation}
Every collection of representations as in \eqref{string} is called \emph{a string}. Each string is actually produced by induction from a scattered representation of $L(s)$ tensored with unitary characters.

Therefore, to completely determine the set $\widehat{E_7}^d$, it suffices to pin down the scattered representations of complex $E_7$ and its proper Levi subgroups. An algorithm for computing the scattered representations of $G$ is provided in \cite[Section 4]{D2}. It will be improved for complex $E_7$ in Section \ref{sec-E7-scattered}.

Our main result is the following.

\medskip
\noindent\textbf{Theorem A.}
\emph{The set $\widehat{E}_7^{\mathrm{d}}$ consists of $66$ scattered representations (see Tables \ref{table-E7-scattered-part-1} and \ref{table-E7-scattered-part-2}) whose spin-lowest $K$-types are all unitarily small, and $578$ strings of representations (see Section \ref{sec-E7-string}).
Moreover, each representation $\pi\in\widehat{E}_7^{\mathrm{d}}$ has a unique spin-lowest $K$-type occurring with multiplicity one.}
\medskip

In Theorem A, the notion unitarily small (\emph{u-small} for short henceforth) $K$-type was introduced by Salamanca-Riba and Vogan \cite{SV}.
Recall that Conjecture 3.4 of \cite{BP} asserts that any Dirac series of $G$ should have a unique spin lowest $K$-type which occurs exactly once, and Conjecture C of \cite{DD} further asserts that the unique spin lowest $K$-type of any scattered representation of $G$ should be u-small. By Theorem A and \cite[Lemma 3.4]{DW1}, both of them now hold for complex $E_7$.

We mention one more corollary of Theorem A.

\begin{cor}\label{cor-main}
Conjecture 1.1 of \cite{BP} holds for complex $E_7$. That is, any Dirac series of complex $E_7$ is either a unipotent representation, or parabolically induced from a unipotent representation  with non-zero Dirac cohomology tensored with unitary characters.
\end{cor}

Unipotent representations of $G$ are conjecturally the building blocks of the unitary dual of $G$. Under the perspective of \cite{B2}, we will give a formulation of unipotent representation relevant to this paper in Section \ref{sec-unip}, and see how it is related to the work of Barbasch and Vogan \cite{BV}, McGovern \cite{Mc}, along with the recent work of Losev, Mason-Brown and Matvieievskyi \cite{LMM}, \cite{MM}.

\medskip
The paper is organized as follows: Section 2 collects necessary preliminaries. Section 3 (resp., Section 4) descries the string part (resp., the scattered part) of $\widehat{E}_7^{\mathrm{d}}$. We realize 63 scattered members in $\widehat{E}_7^{\mathrm{d}}$ as string limits in Section 5, thus providing interesting links between the string part and the scattered part. The three scattered members which can \emph{not} be realized as string limits are unipotent representations (not necessarily special). They are discussed in Section 6. The last section aims to develop effective ways to  determine the spin LKTs. For string limits, the branching result from \cite{BDW} and a technique suggested by Vogan are helpful.

Among all the exceptional Lie groups, complex $E_8$ is the hardest one from the computational point of view.
It is conceivable that our journey here will help people to classify the Dirac series of complex $E_8$ in the future.

\section{Preliminaries}
Let us adopt the setting of the introduction, and assume that $G$ is a connected complex simple Lie group (viewed as a real Lie group).

\subsection{Zhelobenko classification}\label{sec-Zh}
Fix a Borel subgroup $B$ of $G$
containing $H$. Put
$$
\Delta^{+}(\frg_0, \frh_0)=\Delta(\frb_0, \frh_0).
$$
Let $\rho$ be the half sum of roots in $\Delta^{+}(\frg_0, \frh_0)$. Set
$$
\Delta^+(\frg, \frh)=\Delta^+(\frg_0, \frh_0) \times \{0\} \cup \{0\} \times (-\Delta^+(\frg_0, \frh_0)).$$
Let $\rho_{\frg}$ be the half sum of roots in $\Delta^+(\frg, \frh)$.
The restrictions to $\frt\cong \{(x, -x)\mid x\in \frh_0\}$ of these positive roots are
$$
\Delta^+(\frg, \frt)=\Delta^+(\frk, \frt)\cup\Delta^+(\frp, \frt).
$$
Denote by $\rho_c$ (resp., $\rho_n$) the half-sum of roots in $\Delta^+(\frk, \frt)$ (resp., $\Delta^+(\frp, \frt)$). Using the identifications in \eqref{identifications}, we have that
\begin{equation}\label{half-sums}
\rho_{\frg}=(\rho, -\rho), \quad \rho_c=\rho_n=\rho.
\end{equation}
We may and we will identify a $K$-type ($\widetilde{K}$-type, $\frk$-type, etc) with its highest weight, which is necessarily dominant integral with respect to $\Delta^{+}(\frk, \frt)$.

Let $(\lambda_{L}, \lambda_{R})\in \frh_0^{*}\times
\frh_0^{*}$ be such that $\lambda_{L}-\lambda_{R}$ is
a weight of a finite dimensional holomorphic representation of $G$.
 We view $(\lambda_L, \lambda_R)$ as a real-linear functional on $\frh$ by \eqref{identifications}, and write $\bbC_{(\lambda_L, \lambda_R)}$ as the character of $H$ with differential $(\lambda_L, \lambda_R)$ (which exists). By \eqref{identifications} again, we have
$$
\bbC_{(\lambda_L, \lambda_R)}|_{T}=\bbC_{\lambda_L-\lambda_R}, \quad \bbC_{(\lambda_L, \lambda_R)}|_{A}=\bbC_{\lambda_L+\lambda_R}.
$$
Extend $\bbC_{(\lambda_L, \lambda_R)}$ to a character of $B$, and put $$X(\lambda_{L}, \lambda_{R})
:=K\mbox{-finite part of Ind}_{B}^{G}
(
\bbC_{(\lambda_L, \lambda_R)}
) \mbox{ (normalized induction)}.
$$

Using Frobenius reciprocity, the $K$-type with extremal weight $\lambda_{L}-\lambda_{R}$
occurs with multiplicity one in
$X(\lambda_{L}, \lambda_{R})$. Let
$J(\lambda_L,\lambda_R)$ be the unique subquotient of
$X(\lambda_{L}, \lambda_{R})$ containing the
$K$-type $\{\lambda_L -\lambda_R\}$ (recall that this notation stands for the unique dominant weight to which $\lambda_L-\lambda_R$ is conjugate under the Weyl group action).

\begin{thm}\label{thm-Zh} {\rm (Zhelobenko \cite{Zh})}
In the above setting, we have that
\begin{itemize}
\item[a)] Every irreducible admissible ($\frg$, $K$)-module is of the form $J(\lambda_L,\lambda_R)$.
\item[b)] Two such modules $J(\lambda_L,\lambda_R)$ and
$J(\lambda_L^{\prime},\lambda_R^{\prime})$ are equivalent if and
only if there exists $w\in W$ such that
$w\lambda_L=\lambda_L^{\prime}$ and $w\lambda_R=\lambda_R^{\prime}$.
\item[c)] $J(\lambda_L, \lambda_R)$ admits a nondegenerate Hermitian form if and only if there exists
$w\in W$ such that $w(\lambda_L-\lambda_R) =\lambda_L-\lambda_R , w(\lambda_L+\lambda_R) = -\overline{(\lambda_L+\lambda_R)}$.
\item[d)] The representation $X(\lambda_{L}, \lambda_{R})$ is tempered if and only if $\lambda_{L}+\lambda_{R}\in i\frh_0^*$. In this case,
$X(\lambda_{L}, \lambda_{R})=J(\lambda_{L}, \lambda_{R})$.
\end{itemize}
\end{thm}

Note that $J(\lambda_L,\lambda_R)$ has lowest $K$-type $\{\lambda_L-\lambda_R\}$ and infinitesimal character the $W\times W$ orbit of $(\lambda_L, \lambda_R)$.

\subsection{Dirac cohomology}\label{sec-Dirac-coho}

 Let $Z_1, \dots, Z_n$ be an orthonormal basis of $\frp_0$ with respect to the inner product induced by the Killing form $\langle \cdot, \cdot \rangle$. The Dirac operator is defined as
\begin{equation}\label{D-square}
D=\sum_{i=1}^{n} Z_i \otimes Z_i  \in U(\frg)\otimes C(\frp).
\end{equation}
It is easy to check that $D$ does not depend on the choice of the orthonormal basis $\{Z_i\}_{i=1}^{n}$ and it is $K$-invariant for the diagonal action of $K$ given by adjoint actions on both factors. By construction, $D^2$ is a natural Laplacian on $G$. Indeed,
\begin{equation}
D^2=\Omega_{\frk_{\Delta}} -\Omega_{\frg}\otimes 1 +
(\|\rho_c\|^2-\|\rho_{\frg}\|^2) 1 \otimes 1.
\end{equation}
Here $\Omega_{\frg}$ (resp., $\Omega_{\frk}$) is the Casimir operator for $\frg$ (resp., $\frk$), and $\Omega_{\frk_{\Delta}}$ is the image of $\Omega_{\frk}$ under a diagonal embedding $\Delta$ of $\frk$ into $U(\frg)\otimes C(\frp)$. See Section 3.1 of \cite{HP2} for more.

Let $\gamma$ be the highest weight of any $\widetilde{K}$-type of $\pi\otimes S$, where $\pi$ is an irreducible \emph{unitary} $(\frg, K)$ module  which has infinitesimal character $\Lambda$. Then the eigenvalue of $D^2$ on $(X\otimes S)(\gamma)$ is $\|\gamma+\rho_c\|^2-\|\Lambda\|^2$. This leads to
Parthasarathy's Dirac operator inequality
\begin{equation}\label{Dirac-inequality}
\|\gamma+\rho_c\|  \geq \|\Lambda\|.
\end{equation}
As been sharpened in Theorem 3.5.2 of \cite{HP2}, equality  holds in \eqref{Dirac-inequality} if and only if $\gamma+\rho_c$ is conjugate to $\Lambda$ under the action of $W(\frg, \frh)$. Indeed, if \eqref{Dirac-inequality} becomes equality on certain $\widetilde{K}$-type $\gamma$ of $\pi\otimes S$, then
$$(\pi\otimes S)(\gamma)\subseteq \ker D^2= H_D(\pi).
$$
Thus by Theorem \ref{thm-HP}, $\gamma + \rho_c$ is conjugate to $\Lambda$ under the action of $W(\frg, \frh)$.

The following result is the Vogan conjecture proven by Huang and Pand\v zi\'c \cite{HP1}. It is foundational for the study of Dirac cohomology. For instance, although $\pi$ is usually infinite-dimensional, it allows one to focus on finitely many $K$-types to compute $H_D(\pi)$.

\begin{thm}{\rm (Huang and Pand\v zi\'c)}\label{thm-HP}
Let $\pi$ be an irreducible ($\frg$, $K$) module.
Assume that the Dirac
cohomology of $\pi$ is nonzero, and that it contains the $\widetilde{K}$-type $E_{\gamma}$ with highest weight $\gamma\in\frt^{*}\subset\frh^{*}$. Then the infinitesimal character $\Lambda$ of $\pi$ is conjugate to
$\gamma+\rho_{c}$ under $W(\frg,\frh)$.
\end{thm}

\subsection{Spin norm and spin-lowest $K$-type}

Inspired by the ideas of Parthasarathy, Vogan, Huang and Pand\v zi\'c, the first named author introduced spin norm and spin-lowest $K$-type in his 2011 HKUST thesis (cf. \cite{D1}). Given a $K$-type $\mu$, its \emph{spin norm} is defined as
\begin{equation}\label{spin-norm}
\|\mu\|_{\rm spin}:=  \|\{\mu-\rho_n\}+\rho_c\|= \|\{\mu-\rho\}+\rho\|.
\end{equation}
The last step uses \eqref{half-sums}. The spin norm $\|\pi\|_{\rm spin}$ of an irreducible $(\frg, K)$ module $\pi$ is defined as the minimum of the spin norm of all the $K$-types showing up in $\pi$. We call $\mu$ a \emph{spin-lowest $K$-type} of $\pi$ if $\mu$ occurs in $\pi$ and $\|\mu\|_{\rm spin}=\|\pi\|_{\rm spin}$.

Now assume that $\pi$ is unitary and has infinitesimal character $\Lambda$. Let $\mu$ be any $K$-type of $\pi$. Note that the spin module $S$ is some copies of the $\frk$-type $\rho$ (see Lemma 2.2 of \cite{BP}). Thus using the knowledge of PRV-component \cite{PRV}, one sees that the Dirac inequality \eqref{Dirac-inequality} can be reformulated as
\begin{equation}\label{Dirac-inequality-reformulated}
\|\mu\|_{\rm spin}  \geq \|\Lambda\|.
\end{equation}
Moreover, results from the previous section guarantee that:
\begin{itemize}
\item[$\bullet$] $\|\pi\|_{\rm spin}\geq \|\Lambda\|$, and equality holds if and only if $H_D(\pi)$ is non-zero;

\item[$\bullet$] $\|\mu\|_{\rm spin}\geq \|\Lambda\|$, and equality holds if and only if $\mu$ contributes to $H_D(\pi)$;

\item[$\bullet$] If $H_D(\pi)\neq 0$, it is exactly the spin-lowest $K$-types of $\pi$ that contribute to $H_D(\pi)$.
\end{itemize}
These properties motivate us to pay close attention to spin norm and spin-lowest $K$-types when classifying Dirac series.

\subsection{Vogan pencil}
Let $\beta$ be the highest root of $\Delta^{+}(\frg_0, \frh_0)$.
As a special case of Lemma 3.4 and Corollary 3.5 of \cite{Vog81}, for any infinite-dimensional irreducible ($\frg$, $K$)-module $\pi$, there is a unique set
$$\{\mu_{\fri} \,|\, \fri \in  I \} \subseteq
i \frt_{0}^{*}$$ of  dominant integral weights such that all the
$K$-types of $\pi$ are precisely
$$\{\mu_{\fri} +n \beta
\,|\, \fri \in  I, n\in\bbN \}.$$
We call a set of $K$-types
\begin{equation}\label{P-delta}
P(\delta):=\{\delta +n \beta \,|\, n\in\bbN \}
\end{equation}
a \emph{Vogan pencil}. For instance, $P(0)$ denotes the pencil starting from the trivial $K$-type. Put
 \begin{equation}\label{P-mu-prime}
P_{\delta}:=\min \{\|\delta +n \beta \|_{\mathrm{spin}}\,|\, n\in\bbN \}.
\end{equation}
That is, the minimum spin norm of the $K$-types on $P(\delta)$ is denoted as $P_{\delta}$.

\section{The string part of $\widehat{E}_7^{\mathrm{d}}$}\label{sec-E7-string}
To understand the string part of $\widehat{E}_7^{\mathrm{d}}$, we should obtain the scattered parts of $\widehat{A}_i^{\mathrm{d}}$ for $1\leq i\leq 6$ and $\widehat{D}_j^{\mathrm{d}}$ for $4\leq j\leq 6$. Here by $D_j$ we actually mean the group $Spin(2j, \bbC)$. Let us present  the scattered parts of $\widehat{A}_6^{\mathrm{d}}$ and $\widehat{D}_6^{\mathrm{d}}$ in
Tables \ref{table-A6-scattered-part} and \ref{table-D6-scattered-part}, respectively.
Here we present the \emph{folded version} as in \cite[Section 3]{D2}:
If the longest element $w_0$ of $W$ is not equal to $-1$, then it can happen that  $J(\lambda, -s\lambda)$ and $J(-w_0\lambda, w_0 s\lambda)$ are different modules, while both of them are Dirac series. In such a case, we will only present one of them and put a star on its parameter since the other one can be recovered without much difficulty.

Therefore, in both
Table \ref{table-A6-scattered-part} and Table \ref{table-D6-scattered-part}, the number of scattered representations is the number of representations in the table plus the number of stars in the table.
It follows that $N_{A_6}=32$ and $N_{D_6}=34$. By Lemma 3.4 of \cite{DW1}, we conclude that each spin-lowest $K$-type of any scattered representation in Table \ref{table-A6-scattered-part} or Table \ref{table-D6-scattered-part} is u-small. Note that scattered representations of $\widehat{A_i}^d$ for $1\leq i\leq 5$ and $\widehat{D_j}^d$ for $4\leq j\leq 5$ have been recorded in \cite{D2}.

We remark that the scattered representations of complex classical Lie groups have been studied in \cite{DW1} and \cite{DW2}. Be careful that the group of type $D_j$ in \cite{DW2} actually means $SO(2j, \bbC)$. In particular, a representation in Table \ref{table-D6-scattered-part} passes to a representation of $SO(12, \bbC)$ if and only if $\lambda_5=\lambda_6$ holds there; otherwise, it is a genuine representation of $Spin(12, \bbC)$.  Genuine unitary representations of
$Spin(n, \bbC)$ with non-zero Dirac cohomology will be classified in a subsequent work.

\begin{table}
\centering
\caption{The scattered part of $\widehat{A}_6^{\mathrm{d}}$ (folded version)}
\begin{tabular}{l|c|c|c|r}
$s\rho$ &   $\lambda$   & spin LKT & mult \\
\hline
$[3,-5,3,3,-5,3]$ & $\frac{\rho}{2}$ & $\rho$  & $1$ \\
$[-3, 1, 4, -6, 5, -2]^*$ & $[\frac{1}{2}, \frac{1}{2}, \frac{1}{2}, \frac{1}{2}, \frac{1}{2}, 1]$ & $[1, 1, 1, 2, 0, 1]$& $1$  \\
$[-2,-1,5,-6,5,-2]^*$ & $[1, 1, \frac{1}{2}, \frac{1}{2}, \frac{1}{2}, 1]$ & $[1, 0, 0, 4, 0, 1]$ & $1$  \\
$[-4,2,-1,4,-6,4]^*$ & $[\frac{1}{2}, \frac{1}{2}, 1, \frac{1}{2}, \frac{1}{2}, \frac{1}{2}]$ & $[2, 1, 0, 1, 2, 1]$ & $1$ \\
$[-2,-1,-1,5,-6,4]^*$ & $[1, \frac{1}{2}, \frac{1}{2}, \frac{1}{2}, \frac{1}{2}, \frac{1}{2}]$ & $[1, 1, 0, 2, 1, 1]$& $1$ \\
$[-2,-1,-1,5,-6,4]^*$ & $[1, 1, 1, \frac{1}{2}, \frac{1}{2}, \frac{1}{2}]$ & $[1, 0, 0, 0, 4, 1]$& $1$  \\
$[-5,3,-1,-1,3,-5]$ & $\frac{\rho}{2}$ & $\rho$ & $1$  \\
 $[-5,3,-1,-1,3,-5]$ & $[\frac{1}{2}, \frac{1}{2}, 1, 1, \frac{1}{2}, \frac{1}{2}]$ & $[3,1,0,0,1,3]$ & $1$ \\
$[-3,1,-3,1,3,-5]^*$ & $\frac{\rho}{2}$ & $\rho$ & $1$ \\
$[-2,3,-5,3,-4,2]^*$ & $[1, \frac{1}{2}, \frac{1}{2}, \frac{1}{2}, \frac{1}{2}, \frac{1}{2}]$ & $[1, 1, 0, 2, 1, 1]$ & $1$ \\
$[-2,-2,1,-2,4,-5]^*$ & $[1, \frac{1}{2}, \frac{1}{2}, 1, \frac{1}{2}, \frac{1}{2}]$ & $[2, 1, 0, 0, 2, 2]$ & $1$\\
$[-2,-2,1,-2,4,-5]^*$ & $[1, 1, \frac{1}{2}, \frac{1}{2}, \frac{1}{2}, \frac{1}{2}]$ & $[1, 1, 0, 0, 3, 1]$ & $1$ \\
$[-2,-1,-1,-1,4,-5]^*$ & $[1, 1, 1, 1, \frac{1}{2}, \frac{1}{2}]$ & $[0,0,0,0,0,4]$ & $1$ \\
$[-2,-2,3,-4,2,-3]^*$ & $[1, 1, \frac{1}{2}, \frac{1}{2}, \frac{1}{2}, \frac{1}{2}]$ & $[1, 1, 0, 0, 3, 1]$ & $1$\\
$[-1,-3,1,1,-3,-1]$ & $[1, \frac{1}{2}, \frac{1}{2}, \frac{1}{2}, \frac{1}{2}, 1]$ & $[1, 2, 0, 0, 2, 1]$& $1$\\
$[-1,-2,-1,2,-3,-1]^*$ & $[1, 1, 1, \frac{1}{2}, \frac{1}{2}, 1]$ & $[2, 0, 0, 0, 0, 4]$ & $1$ \\
$-\rho$ & $[1, \frac{1}{2}, \frac{1}{2}, \frac{1}{2}, \frac{1}{2}, 1]$ & $[1, 2, 0, 0, 2, 1]$ & $1$ \\
$-\rho$ & $[1, 1, \frac{1}{2}, \frac{1}{2}, 1, 1]$ & $[3, 0, 0, 0, 0, 3]$& $1$ \\
$-\rho$ & $\frac{\rho}{2}$ & $\rho$ & $1$\\
$-\rho$ & $\rho$ & $[0,0,0,0,0,0]$ & $1$
\end{tabular}
\label{table-A6-scattered-part}
\end{table}

\begin{table}
\centering
\caption{The scattered part of $\widehat{D}_6^{\mathrm{d}}$ (folded version)}
\begin{tabular}{l|c|c|c|r}
$s\rho$ &   $\lambda$   & spin LKT & mult \\
\hline
$[-3, 1, 5, -7, 5, 5]$ & $\frac{\rho}{2}$ & $\rho$ & $1$\\
$[-2, -1, 6, -7, 5, 5]$ & $[1,1,\frac{1}{2},\frac{1}{2},\frac{1}{2},\frac{1}{2}]$ & $[1, 0, 0, 3, 1, 1]$ & $1$ \\
$[-2, 7, -8, 6, -1, -1]$ & $[1,\frac{1}{2},\frac{1}{2},\frac{1}{2},1,1]$ & $[1, 0, 4, 1, 0, 0]$ & $1$ \\
$[4, -6, 5, 4, -6, -2]^*$ & $[\frac{1}{2},\frac{1}{2},\frac{1}{2},\frac{1}{2},\frac{1}{2},1]$ & $[1, 2, 0, 1, 2, 1]$ & $1$ \\
$[7, -9, 7, -3, 1, 1]$ & $\frac{\rho}{2}$ & $\rho$ & $1$ \\
$[-5, 4, -2, 7, -8, -2]^*$ & $[\frac{1}{2},\frac{1}{2},1,\frac{1}{2},\frac{1}{2},1]$ & $[3, 0, 1, 0, 4, 1]$ & $1$ \\
$[5, -7, 5, -7, 5, 5]$ & $\frac{\rho}{2}$ & $\rho$ & $1$ \\
$[7, -9, 7, -1, -1, -1]$ & $[\frac{1}{2},\frac{1}{2},\frac{1}{2},1,1,1]$ & $[1, 5, 1, 0, 0, 0]$ & $1$ \\
$[-1, -3, 1, 7, -9, -1]^*$ & $[1,\frac{1}{2},\frac{1}{2},\frac{1}{2},\frac{1}{2},1]$ & $[2, 0, 1, 0, 3, 2]$ & $1$ \\
$[-6, 4, 5, -6, -2, 4]^*$ & $[\frac{1}{2},\frac{1}{2},\frac{1}{2},\frac{1}{2},1,\frac{1}{2}]$ & $[2, 0, 2, 0, 2, 1]$ & $1$ \\
$[-1, -2, -1, 8, -9, -1]^*$ & $[1,1,1,\frac{1}{2},\frac{1}{2},1]$ & $[0, 1, 0, 0, 7, 0]$ & $1$ \\
$[-9, 7, -5, 3, -1, -1]$ & $[\frac{1}{2},\frac{1}{2},\frac{1}{2},\frac{1}{2},1,1]$ & $[1, 3, 1, 1, 0, 0]$ & $1$ \\
$[-2, -1, 6, -8, -1, 7]^*$ & $[1,1,\frac{1}{2},\frac{1}{2},1,\frac{1}{2}]$ & $[0, 2, 0, 0, 5, 0]$ & $1$\\
$[-7, 5, -7, 5, -1, -1]$ & $[\frac{1}{2},\frac{1}{2},\frac{1}{2},\frac{1}{2},1,1]$ & $[1, 3, 1, 1, 0, 0]$ & $1$ \\
$[-9, 7, -1, -3, 1, 1]$ & $[\frac{1}{2},\frac{1}{2},1,\frac{1}{2},\frac{1}{2},\frac{1}{2}]$ & $[2, 1, 2, 1, 0, 0]$ & $1$ \\
$[-1, -7, 5, -3, 1, 1]$ & $[1,\frac{1}{2},\frac{1}{2},\frac{1}{2},\frac{1}{2},\frac{1}{2}]$ & $[1, 1, 1, 2, 0, 0]$ & $1$ \\
$[-9, 7, -1, -1, -1, -1]$ & $[\frac{1}{2},\frac{1}{2},1,1,1,1]$ & $[7, 1, 0, 0, 0, 0]$ & $1$ \\
$[-1, -5, 3, -5, 3, 3]$ & $[1,\frac{1}{2},\frac{1}{2},\frac{1}{2},\frac{1}{2},\frac{1}{2}]$ & $[1, 1, 1, 2, 0, 0]$ & $1$ \\
$[-2, 4, -5, 4, -6, -2]^*$ & $[1,\frac{1}{2},\frac{1}{2},\frac{1}{2},\frac{1}{2},1]$ & $[0, 2, 0, 1, 0, 3]$ & $1$\\
$[-1, -7, 5, -1, -1, -1]$ & $[1,\frac{1}{2},\frac{1}{2},1,1,1]$ & $[5, 2, 0, 0, 0, 0]$ & $1$ \\
$[-1, -1, -5, 3, -1, -1]$ & $[1,1,\frac{1}{2},\frac{1}{2},1,1]$ & $[3, 3, 0, 0, 0, 0]$ & $1$ \\
$[1, -3, 1, -3, -1, 1]^*$ & $[\frac{1}{2},\frac{1}{2},\frac{1}{2},\frac{1}{2},1,\frac{1}{2}]$ & $[0, 2, 0, 2, 0, 1]$& $1$ \\
$[-1, -1, -1, -3, 1, 1]$ & $\frac{\rho}{2}$ & $\rho$ & $1$ \\
$[-1, -1, -1, -3, 1, 1]$ & $[1,1,1,\frac{1}{2},\frac{1}{2},\frac{1}{2}]$ & $[1, 4, 0, 0, 0, 0]$& $1$ \\
$-\rho$ & $[1,\frac{1}{2},\frac{1}{2},\frac{1}{2},\frac{1}{2},\frac{1}{2}]$ & $[1, 1, 1, 2, 0, 0]$& $1$ \\
$-\rho$ & $\rho$ & [0,0,0,0,0,0]& $1$
\end{tabular}
\label{table-D6-scattered-part}
\end{table}

It is not hard to equip these scattered representations into strings of $\widehat{E}_7^{\mathrm{d}}$. Let us see some examples.

\begin{example}\label{exam-HA}
Consider the extreme case that $s=e$. Then $I(s)=\{1, 2, 3, 4, 5, 6, 7\}$ and $L(s)=H$. Recalling \eqref{Lambda-s}, one computes that
$$
\Lambda(e)=\{\lambda=[\lambda_1, \dots, \lambda_7] \mid \lambda_i\in\frac{1}{2}\mathbb{P}\}.
$$
Each $\lambda\in\Lambda(e)$ gives a representation $J(\lambda, -e\lambda)=J(\lambda, -\lambda)$ which is a Dirac series. Indeed, $J(\lambda, -\lambda)$ has LKT $2\lambda$, which contributes to the Dirac cohomology since
$$
\|2\lambda\|_{\rm spin}=\|\{2\lambda -\rho\} + \rho\|=\|2\lambda\|=\|\Lambda\|.
$$
All these modules form the single string (recall equation \eqref{string}) $(e, [\lambda_1, \dots, \lambda_7])$ with each $\lambda_i$ running over $\frac{1}{2}\mathbb{P}$. The underlying representation of this string is the trivial one of $H$. \hfill\qed
\end{example}
\begin{rmk}
In general, all the tempered representations of $G$ with non-zero Dirac cohomology is equipped into exactly one string.
\end{rmk}

\begin{example}\label{exam-D6}
The first representation in Table \ref{table-D6-scattered-part} has $s=s_1s_2s_4s_3s_2s_ 1s_5s_6s_4$ and
$$
\lambda=[\frac{1}{2}, \frac{1}{2}, \frac{1}{2}, \frac{1}{2}, \frac{1}{2}, \frac{1}{2}].
$$
Note that $\lambda+s\lambda=[-1, 1, 3, -3, 3, 3]$ and $\lambda-s\lambda=[2, 0, -2, 4, -2, -2]$. Recall that $s\rho(D_6)=[-3, 1, 5, -7, 5, 5]$. We \emph{embed} $D_6$ into $E_7$ in the following way: label the simple root $\alpha_i$ in Fig.~\ref{E7-Dynkin} as the $(8-i)$-th simple root of $D_6$ for $2\leq i\leq 7$. Then the counterpart of the involution $s$ in $W(D_6)$ is the following involution $s'$ in $W(E_7)$:
$$
s'=s_7s_6s_4s_5s_6s_7s_3s_2s_4,
$$
and that
$$
s'\rho(E_7)=[\textbf{3}, 5,5,-7,5,1,-3].
$$
Here the non-bolded coordinate  comes from those of $s\rho(D_6)$ via the permutation determined by the embedding of $D_6$ into $E_7$, while the bolded coordinate corresponds to the simple root $\alpha_1$ in Fig.~\ref{E7-Dynkin}.

Now $I(s^\prime)=\{1\}$ and the semisimple factor of $L(s^\prime)$ is of type $D_6$.
By Proposition 3.5 of \cite{DD},  the representation $J(\lambda, -s\lambda)$ of $D_6$ produces the following  string of $\widehat{E}_7^{\mathrm{d}}$ via induction:
$$
J(\lambda^\prime, -s^\prime \lambda^\prime),
$$
where $\lambda^\prime:=[\frac{a}{2}, \frac{1}{2}, \frac{1}{2}, \frac{1}{2}, \frac{1}{2}, \frac{1}{2}, \frac{1}{2}]$ with $a$ running over $\mathbb{P}$.
Note that
$$
\lambda^\prime+ s^\prime \lambda^\prime=[a +1, 3, 3, -3, 3, 1, -1], \quad \lambda^\prime- s^\prime \lambda^\prime=[-1, -2, -2, 4, -2, 0, 2].
$$
\hfill\qed
\end{example}

Embedding $D_6$ into $E_7$ as in Example \ref{exam-D6}, Table \ref{table-D6-scattered-part} produces $N_{D_6}=34$ strings in $\widehat{E}_7^{\mathrm{d}}$. Similarly, we embed $A_6$ into $E_7$ by labelling the simple root $\alpha_i$ in Fig.~\ref{E7-Dynkin} as the $(i-1)$-th simple root of $A_6$ for $3\leq i\leq 7$, and labelling $\alpha_1$ in Fig.~\ref{E7-Dynkin} as the first simple root of $A_6$. Then Table \ref{table-A6-scattered-part} produces $N_{A_6}=32$ strings in $\widehat{E}_7^{\mathrm{d}}$.

We will not explicitly give all the strings of $\widehat{E}_7^{\mathrm{d}}$. Instead, let us count them. For this purpose, let us firstly mention that by Section 3 of \cite{DW1}, for any $n\in\mathbb{P}$, we have that
$N_{A_{n}}=2^{n-1}$. Moreover, we recall that $N_{D_4}=9$, $N_{D_5}=17$ and $N_{D_6}=34$.
Now denote by $N_i$ the number of $I$-strings in $\widehat{E}_7^{\mathrm{d}}$ such that $|I|=7-i$. Then $N_0 =1$ (see Example \ref{exam-HA}), and
\begin{align*}
N_1 &=7 N_{A_1}=7\times 1 =7,\\
N_2 &=6 N_{A_2}+15 N_{A_1} N_{A_1}=6\times 2 + 15 \times 1 \times 1= 27.
\end{align*}
In the following, we will omit the factors $N_{A_1}$ (which equals to $1$) to save some space.
\begin{align*}
N_3 &=6 N_{A_3} +18 N_{A_2} N_{A_1}  + 11 N_{A_1} N_{A_1} N_{A_1}\\
    &=6\times 4 + 18\times 2 + 11=71,\\
N_4 &= N_{D_4} + 5 N_{A_4}+11 N_{A_3} N_{A_1}+12 N_{A_2} N_{A_1} N_{A_1} + 4 N_{A_2} N_{A_2} + 2 N_{A_1} N_{A_1} N_{A_1} N_{A_1}\\
    &=9+5\times 8+ 11 \times 4 +12 \times 2+4 \times 2 \times 2 +2 =135,\\
N_5 &=2 N_{D_5}+3 N_{A_5}+ N_{D_4}N_{A_1}+ 5 N_{A_4} N_{A_1}
+ 3 N_{A_3} N_{A_2}   + 3 N_{A_3} N_{A_1}N_{A_1} + 3 N_{A_2} N_{A_2}N_{A_1}\\
&
+ N_{A_2}N_{A_1}N_{A_1}N_{A_1} \\
&=2\times 17 + 3\times 16 + 9 +5\times 8 + 3\times 4\times 2 + 3\times 4 + 3\times 2\times 2 + 2 =181, \\
N_6 &= N_{D_6}+ N_{A_6}+ N_{A_5}N_{A_1}+ N_{A_3}N_{A_2}N_{A_1}+N_{A_4}N_{A_2}+N_{D_5}N_{A_1}+N_{E_6}\\
    &=34+32+16+4\times 2+ 8\times 2 + 17 + 33=156.
\end{align*}
Here recall that $N_G$ is the cardinality of the scattered part of $\widehat{G}^{\mathrm{d}}$. Therefore, $\widehat{E}_7^{\mathrm{d}}$ contains $\sum_{i=0}^{6} N_i=578$ strings in total.

\section{The scattered part of $\widehat{E}_7^{\mathrm{d}}$}\label{sec-E7-scattered}

Similar to Section 4 of \cite{D2}, we proceed as follows  to sieve out all the \emph{non-trivial} scattered representations of complex $E_7$:
\begin{itemize}
\item[$\bullet$] collect the finitely many $\lambda\in \Lambda(s)$ (recall equation \eqref{Lambda-s}) such that
$\|\lambda-s\lambda\|^2\leq 464$ and that
$\|2\lambda\|^2 \leq P_{\{\lambda+s\lambda\}}^2$ (see \eqref{P-mu-prime}).

\item[$\bullet$] For these $\lambda$ surviving from the above step (if exist), use  \texttt{atlas} \cite{ALTV,At} to study the unitarity and $K$-types of $J(\lambda, -s\lambda)$.
\end{itemize}
In the first step, note that the earlier upper bound for $\|\lambda-s\lambda\|^2$ is
$\|2\rho\|^2$, which equals to $798$ for complex $E_7$. This improvement is guaranteed by Theorem 1.1 of \cite{D3}.

The relevant files are available from the following link:
\begin{verbatim}
https://www.researchgate.net/publication/353914493_E7-Files
\end{verbatim}

\begin{example}\label{exam-E7-first-scattered}
Consider the involution $s=s_1 s_4 s_2 s_3 s_1 s_5 s_6 s_7 s_6 s_5 s_4$. Note that $s\rho(E_7)=[-2,6,7,-8,6,1,-3]$. Carrying out the first step above leaves us with the following $6$ candidate representations:
\begin{align*}
&[\frac{1}{2}, \frac{1}{2}, 1, \frac{1}{2}, \frac{1}{2}, \frac{1}{2}, \frac{1}{2}], \quad
[1, \frac{1}{2}, \frac{1}{2}, \frac{1}{2}, \frac{1}{2}, \frac{1}{2}, \frac{1}{2}], \quad
[1, \frac{1}{2}, \frac{1}{2}, \frac{1}{2}, \frac{1}{2}, \frac{1}{2}, \frac{3}{2}],\\
&[1, \frac{1}{2}, \frac{1}{2}, \frac{1}{2}, 1, \frac{1}{2}, 1], \quad
[1, \frac{3}{2}, \frac{1}{2}, \frac{1}{2}, \frac{1}{2}, \frac{1}{2}, \frac{1}{2}], \quad
[\frac{3}{2}, \frac{1}{2}, 1, \frac{1}{2}, \frac{1}{2}, \frac{1}{2}, \frac{1}{2}].
\end{align*}
Only the $\lambda=[1, \frac{1}{2}, \frac{1}{2}, \frac{1}{2}, \frac{1}{2}, \frac{1}{2}, \frac{1}{2}]$ produces a unitary representation. Note that
\begin{align*}
\lambda + s \lambda=[0, 4, 4, -4, 4, 1, -1],\quad
\lambda - s \lambda=[2, -3, -3, 5, -3, 0, 2].
\end{align*}
It has a unique spin-lowest $K$-type $[1,1,0,2,1,1,1]$ whose spin norm equals to $\|2\lambda\|$. Moreover, it occurs exactly once. This is the 1st entry of Table \ref{table-E7-scattered-part-1}.
\hfill\qed
\end{example}

\begin{example}\label{exam-E7-3424}
Consider the involution
$$
s=s_1s_2s_3s_4s_2s_3s_4s_5s_4s_2s_3s_4s_5s_6s_5s_4s_2s_3s_4s_5s_6s_7s_6s_5s_4s_2
s_3s_1s_4s_5s_6s_7.
$$ Note that $s\rho(E_7)=[-17,-1,15,-1,-1,-1,-1]$. Carrying out the first step above leaves us with $241$ candidate representations.
Only  $\lambda=[\frac{1}{2}, 1, \frac{1}{2}, 1, 1, 1, 1]$ produces a unitary representation. Note that
$$
\lambda + s\lambda=[-14, 0, 14, 0, 0, 0, 0], \quad \lambda + s\lambda=[15, 2, -13, 2, 2, 2, 2].
$$
It has a unique spin-lowest $K$-type $[13,0,1,0,0,0,0]$ whose spin norm equals to $\|2\lambda\|$. Moreover, it occurs exactly once. This is the 1st entry of Table \ref{table-E7-scattered-part-2}.
\hfill\qed
\end{example}

\begin{example}\label{exam-spherical}
Consider the longest involution
$w_0=-1$ in $W(E_7)$. Note that $w_0\rho(E_7)=[-1,-1,-1,-1,-1,-1,-1]$. Carrying out the first step above leaves us with $116$ candidate representations.
Only the $\lambda=[\frac{1}{2}, \frac{1}{2}, \frac{1}{2}, \frac{1}{2}, \frac{1}{2}, \frac{1}{2}, 1]$ produces a unitary representation. It has a unique spin-lowest $K$-type $[1,0,1,2,0,2,0]$ whose spin norm equals to $\|2\lambda\|$. Moreover, it occurs exactly once. This is the 21st entry of Table \ref{table-E7-scattered-part-2}.
\hfill\qed
\end{example}

\section{String limits}

This section aims to explain the last columns of Tables \ref{table-E7-scattered-part-1} and \ref{table-E7-scattered-part-2}.

\begin{example}\label{ex-string-limit-E6-first}
The first row of \cite[Table 6]{D2} gives a scattered member of $\widehat{E}_6^{\mathrm{d}}$ with $s=s_4 s_5 s_6 s_5 s_1 s_3 s_2 s_4 s_1$ and $\lambda=[1, \frac{1}{2}, \frac{1}{2}, \frac{1}{2}, \frac{1}{2}, 1]$. Note that $s\rho(E_6)=[-2,5,6,-7, 6, -2]$. This scattered representation in $\widehat{E}_6^{\mathrm{d}}$ produces the following string in $\widehat{E}_7^{\mathrm{d}}$:
\begin{equation}\label{string-E6-first}
J(\lambda^\prime, -s\lambda^\prime), \quad \lambda^\prime=[1, \frac{1}{2}, \frac{1}{2}, \frac{1}{2}, \frac{1}{2}, 1, \frac{g}{2}],
\end{equation}
where $g\in\mathbb{P}$. Note that
$$
\lambda^\prime + s \lambda^\prime=[0, 4, 4, -4, 4, 0, g+2], \quad  \lambda^\prime - s \lambda^\prime=[2, -3, -3, 5, -3, 2, -2].
$$

If we set $g=-1$, the corresponding representation is no longer a member of the string. However, $s_7 ([1, \frac{1}{2}, \frac{1}{2}, \frac{1}{2}, \frac{1}{2}, 1, -\frac{1}{2}])=[1, \frac{1}{2}, \frac{1}{2}, \frac{1}{2}, \frac{1}{2}, \frac{1}{2}, \frac{1}{2}]$; moreover, if we set $s^\prime=s_7 s s_7^{-1}$, then $s^\prime$ is an involution of $W(E_7)$ such that $I(s^\prime)$ is empty and that $s^\prime \rho(E_7)=[-2, 6, 7, -8, 6, 1, -3]$. One sees that $s^\prime$ and $[1, \frac{1}{2}, \frac{1}{2}, \frac{1}{2}, \frac{1}{2}, \frac{1}{2}, \frac{1}{2}]$ give the first row of Table \ref{table-E7-scattered-part-1}. Therefore, we may view the first entry of Table \ref{table-E7-scattered-part-1} (see Example \ref{exam-E7-first-scattered}), a scattered member of $\widehat{E}_7^{\mathrm{d}}$, as a limit of the above string $J(\lambda^\prime, -s\lambda^\prime)$, whose underlying representation is a scattered member of $\widehat{E}_6^{\mathrm{d}}$.

Let us see how to realize the above things in \texttt{atlas}.
\begin{verbatim}
set G=complex(simply_connected (E7))
set x=KGB(G,2903039)
set q=parameter(x,[0,4,4,-4,4,1,-1,  0,0,0,0,0,0,0],
                  [2,-3,-3,5,-3,0,2, 0,0,0,0,0,0,0])
set p=monomials(finalize(q))[0]
p
Value: final parameter(x=17376,lambda=[2,-2,-2,4,-2,0,2,2,0,0,2,0,1,1]/1,
                                   nu=[2,-3,-3,5,-3,0,2,2,-3,-3,5,-3,0,2]/2)
\end{verbatim}
The representation \texttt{p} is the first entry of Table \ref{table-E7-scattered-part-1}.
\begin{verbatim}
set q1=parameter(x,[0,4,4,-4,4,0,3, 0,0,0,0,0,0,0],
                   [2,-3,-3,5,-3,2,-2, 0,0,0,0,0,0,0])
set p1=monomials(finalize(q1))[0]
p1
Value: final parameter(x=6951,lambda=[2,-1,-1,3,-2,2,-1,2,0,0,2,0,1,1]/1,
                                  nu=[2,-3,-3,5,-3,2,-2,2,-3,-3,5,-3,2,-2]/2)
\end{verbatim}
The representation \texttt{p1} is the string member with $g=1$.

\begin{verbatim}
set (Q,pp1)=reduce_good_range(p1)
Q
Value: ([0,1,2,3,4,5,7,8,9,10,11,12],KGB element #6951)
pp1
Value: final parameter(x=2670,lambda=[2,-1,-1,3,-2,2,-10,2,0,0,2,0,1,-8]/1,
                                  nu=[2,-3,-3,5,-3,2,-2,2,-3,-3,5,-3,2,-2]/2)
\end{verbatim}
Let us move the inducing module from $g=1$ to $g=-1$, and then do cohomological induction.
\begin{verbatim}
set ppm1=parameter(x(pp1),lambda(pp1)-[0,0,0,0,0,0,2,0,0,0,0,0,0,0],nu(pp1))
goodness(ppm1,G)
Value: "Fair"
set ind=theta_induce_irreducible(ppm1,G)
ind
Value:
1*parameter(x=17376,lambda=[2,-2,-2,4,-2,0,2,2,0,0,2,0,1,1]/1,
                        nu=[2,-3,-3,5,-3,0,2,2,-3,-3,5,-3,0,2]/2) [411]
monomials(ind)[0]=p
Value: true
\end{verbatim}
Thus the induced module is irreducible, and it is exactly the representation \texttt{p}. \hfill\qed
\end{example}

For convenience, we simply say that the first entry of Table \ref{table-E7-scattered-part-1} is a \textbf{string limit} coming from the first scattered representation of $E_6$ with $g=-1$.

\begin{example}\label{ex-string-limit-E6-last}
The last row of \cite[Table 6]{D2} is the trivial representation of complex $E_6$. Namely, we have that  $s\rho(E_6)=[-1,-1,-1,-1,-1,-1]$ and $\lambda=[1, 1, 1, 1, 1, 1]$. It produces the following string in $\widehat{E}_7^{\mathrm{d}}$:
\begin{equation}\label{string-E6-last}
J(\lambda^\prime, -s\lambda^\prime), \quad \lambda^\prime=[1, 1, 1, 1, 1, 1, \frac{g}{2}],
\end{equation}
where $g\in\mathbb{P}$. Note that
$$
\lambda^\prime + s \lambda^\prime=[0, 0, 0, 0, 0, 0,  g+16], \quad  \lambda^\prime - s \lambda^\prime=[2, 2, 2, 2, 2, 2, -16].
$$
This string produces the following eight limits in total:
\begin{itemize}
\item[$\bullet$] $g=-1$, the 8th entry of Table \ref{table-E7-scattered-part-2};

\item[$\bullet$] $g=-3$, the 9th entry of Table \ref{table-E7-scattered-part-2};

\item[$\bullet$] $g=-5$, the 11th entry of Table \ref{table-E7-scattered-part-2};

\item[$\bullet$] $g=-7$, the 13th entry of Table \ref{table-E7-scattered-part-2};

\item[$\bullet$] $g=-9$, the 14th entry of Table \ref{table-E7-scattered-part-2};

\item[$\bullet$] $g=-11$, the 16th entry of Table \ref{table-E7-scattered-part-2};

\item[$\bullet$] $g=-13$, the 17th entry of Table \ref{table-E7-scattered-part-2};

\item[$\bullet$] $g=-15$, the 20th entry of Table \ref{table-E7-scattered-part-2}.
\end{itemize}
Note that in each case, the inducing module is in the fair range.
\hfill\qed
\end{example}

As recorded in the last columns of Tables \ref{table-E7-scattered-part-1} and \ref{table-E7-scattered-part-2}, we \emph{exhaust} all the string limits coming from the scattered representations of the $E_6$ factor, and the inducing module is always in the fair range. Similar things have been done for the $A_6$ and $D_6$ factor. Let us give a few more examples.

\begin{example}\label{ex-string-limit-A6}
The 13th row of Table \ref{table-A6-scattered-part} gives a scattered member of $\widehat{A}_6^{\mathrm{d}}$ with $s\rho(A_6)=[-2,-1,-1, -1,4,-5]$ and $\lambda=[1, 1, 1, 1,\frac{1}{2}, \frac{1}{2}]$. It produces the following string in $\widehat{E}_7^{\mathrm{d}}$:
\begin{equation}\label{string-A6-13th}
J(\lambda^\prime, -s\lambda^\prime), \quad \lambda^\prime=[1, \frac{b}{2}, 1, 1, 1, \frac{1}{2}, \frac{1}{2}],
\end{equation}
where $b\in\mathbb{P}$. Note that
$$
\lambda^\prime + s \lambda^\prime=[0, b+9, 0, 0, 0, 4, -4], \quad  \lambda^\prime - s \lambda^\prime=[2, -9, 2, 2, 2, -3, 5].
$$
This string gives the following three limits in total:
\begin{itemize}
\item[(a)] $b=-1$, the 17th entry of Table \ref{table-E7-scattered-part-1};

\item[(b)] $b=-3$, the 24th entry of Table \ref{table-E7-scattered-part-1};

\item[(c)] $b=-20$, the 18th entry of Table \ref{table-E7-scattered-part-1}.
\end{itemize}
Note that both (a) and (b) are also string limits coming from the $E_6$ factor, while (c) is not. Moreover, the inducing module for (c) is not in the fair range.

The 13th row of Table \ref{table-A6-scattered-part} is marked with a star, meaning that it has a dual representation. Indeed, it is a scattered member of $\widehat{A}_6^d$ such that $s^*\rho(A_6)=[-5,4,-1,-1,-1,-2]$ and $\lambda^*=[\frac{1}{2}, \frac{1}{2}, 1, 1, 1, 1]$. It produces the following string in $\widehat{E}_7^{\mathrm{d}}$:
\begin{equation}\label{string-A6-13th-dual}
J(\lambda^{\prime\prime}, -s^*\lambda^{\prime\prime}), \quad \lambda^{\prime\prime}=[\frac{1}{2}, \frac{b}{2}, \frac{1}{2}, 1, 1, 1,  1],
\end{equation}
where $b\in\mathbb{P}$.
Note that
$$
\lambda^{\prime\prime} + s^* \lambda^{\prime\prime}=[-4, b+8, 4, 0, 0, 0, 0], \quad  \lambda^{\prime\prime} - s^* \lambda^{\prime\prime}=[5, -8, -3, 2, 2, 2, 2].
$$
This string gives the following three limits in total:
\begin{itemize}
\item[(a')] $b=-1$, the representation (c);

\item[(b')] $b=-18$, the representation (b);

\item[(c')] $b=-20$, the representation (a).
\end{itemize}
Note that the inducing module for (a') is now in the fair range.  \hfill\qed
\end{example}

\begin{example}\label{ex-string-limit-D6}
The 19th row of Table \ref{table-D6-scattered-part} gives a scattered member of $\widehat{D}_6^{\mathrm{d}}$ with $s\rho(D_6)=[-2,4,-5,4,-6,-2]$ and $\lambda=[1, \frac{1}{2}, \frac{1}{2}, \frac{1}{2}, \frac{1}{2}, 1]$. It produces the following string in $\widehat{E}_7^{\mathrm{d}}$:
\begin{equation}\label{string-D6-19th}
J(\lambda^\prime, -s\lambda^\prime), \quad \lambda^\prime=[\frac{a}{2}, 1, \frac{1}{2}, \frac{1}{2}, \frac{1}{2}, \frac{1}{2}, 1],
\end{equation}
where $a\in\mathbb{P}$. This string gives no limit.

On the other hand, the dual representation has $s^*\rho(D_6)=[-2,4,-5,4,-2,-6]$ and $\lambda^*=[1, \frac{1}{2}, \frac{1}{2}, \frac{1}{2}, 1, \frac{1}{2}]$. It produces the following string in $\widehat{E}_7^{\mathrm{d}}$:
\begin{equation}\label{string-D6-19th-dual}
J(\lambda^{\prime\prime}, -s^*\lambda^{\prime\prime}), \quad \lambda^{\prime\prime}=[\frac{a}{2}, \frac{1}{2}, 1,  \frac{1}{2}, \frac{1}{2}, \frac{1}{2}, 1],
\end{equation}
where $a\in\mathbb{P}$.
Note that
$$
\lambda^{\prime\prime} + s^* \lambda^{\prime\prime}=[a+8, -3,0,3,-3,3,0], \quad  \lambda^{\prime\prime} - s^* \lambda^{\prime\prime}=[-8,4,2,-2,4,-2,2].
$$
When $a=-1$, this string gives the 31st entry of of Table \ref{table-E7-scattered-part-1}. \hfill\qed
\end{example}

\begin{example}\label{ex-string-limit-D5A1}
Consider the maximal Levi obtained by removing the simple root $\alpha_6$. Put trivial representations on both the $D_5$ factor and the $A_1$ factor. Correspondingly, we have the string $J(\lambda, -s\lambda)$ in $\widehat{E}_7^{\mathrm{d}}$ such that
\begin{equation}
\lambda+s\lambda=[0, 0, 0, 0, 0,  2 f + 11, 0], \quad \lambda+s\lambda=[2, 2, 2, 2, 2, -11, 2]
\end{equation}
where $f\in\mathbb{P}$.
This string gives the following three limits in total:
\begin{itemize}
\item[(a)] $f=-1$, the 28th entry of Table \ref{table-E7-scattered-part-1};

\item[(b)] $f=-3$, the 36th entry of Table \ref{table-E7-scattered-part-1};

\item[(c)] $f=-5$, the last entry of Table \ref{table-E7-scattered-part-1}.
\end{itemize}

We may also put the scattered representation of $\widehat{D}_5^d$ sitting in the second row of \cite[Table 5]{D2}, and put the trivial representation on the $A_1$ factor. They produce the string $J(\lambda', -s'\lambda')$ in $\widehat{E}_7^{\mathrm{d}}$ such that
\begin{equation}
\lambda'+s'\lambda'=[4, 0, -4, 4, 0, 2 f +4 , 0], \quad \lambda'+s'\lambda'=[-3, 2, 5, -3, 2, -4, 2].
\end{equation}
where $f\in\mathbb{P}$. Setting $f=-1$ gives the third entry of Table \ref{table-E7-scattered-part-1}.
 \hfill\qed
\end{example}

\medskip

\begin{table}
\centering
\caption{The scattered representations in $\widehat{E}_7^{\mathrm{d}}$ (part one)}
\begin{tabular}{r|c|c|c|c}
$s\rho$ &   $\lambda$   & spin LKT & mult & string limits \\
\hline
$[-2,6,7,-8,6,1,-3]$ & $[1, \frac{1}{2}, \frac{1}{2}, \frac{1}{2}, \frac{1}{2}, \frac{1}{2}, \frac{1}{2}]$  & $[1, 1, 0, 2, 1, 1, 1]$ & $1$ & $E_6$, 1st, $g=-1$\\

$[-2,6,7,-8,7,-1,-2 ]$ & $[1, \frac{1}{2}, \frac{1}{2}, \frac{1}{2}, \frac{1}{2}, 1, 1]$ & $[1,1,0,4,0,0,1]$ & $1$ &$A_2.A_1.A_3$, $d=-1$\\

$[6,-2,-8,6,5,-6,4 ]$ & $[\frac{1}{2}, 1, \frac{1}{2}, \frac{1}{2}, \frac{1}{2}, \frac{1}{2}, \frac{1}{2}]$ & $[1,1,2,1,0,2,1 ]$ & $1$ &  $D_5$ 2nd.$A_1$, $f=-1$\\

$[-5,-7,3,5,5,-7,5 ]$ & $\frac{\rho}{2}$ & $\rho$  & $1$  & $A_6$, 4th, $b=-1$ \\

$[-4,-2,3,7,-9,8,-2 ]$ & $[\frac{1}{2}, 1, \frac{1}{2}, \frac{1}{2}, \frac{1}{2}, \frac{1}{2}, 1 ]$ & $[2,1,0,1,3,0,1 ]$ & $1$  & $D_6$, 3rd*, $a=-1$\\

$[-2,-8,-2,7,6,-7,5 ]$ & $[1, \frac{1}{2}, 1, \frac{1}{2}, \frac{1}{2}, \frac{1}{2}, \frac{1}{2} ]$ & $[1,4,1,0,0,3,1]$ & $1$ & $A_6$, 6th, $b=-1$\\
$[-1,-1,-2,9,-10,9,-2 ]$ & $[1,1,1,\frac{1}{2}, \frac{1}{2}, \frac{1}{2},1]$ & $[0,0,1,0,6,0,1 ]$  & $1$ & $A_4.A_2$, $e=-1$\\

$[8,-2,-10,9,-2,4,-5 ]$ & $[\frac{1}{2}, 1, \frac{1}{2}, \frac{1}{2}, 1, \frac{1}{2}, \frac{1}{2}]$ & $[1,1,4,0,1,0,3]$  & $1$ & $E_6$, 4th*, $g=-1$\\

$[-6,-9,4,8,-2,4,-6 ]$ & $[\frac{1}{2}, \frac{1}{2}, \frac{1}{2}, \frac{1}{2}, 1, \frac{1}{2}, \frac{1}{2} ]$ & $[2,3,1,0,1,1,2 ]$  & $1$ & $A_6$, 8th, $b=-1$\\

$[-8,-1,6,-1,8,-10,8 ]$ & $[\frac{1}{2},1,\frac{1}{2},1,\frac{1}{2},\frac{1}{2},\frac{1}{2}]$ & $[4,0,1,0,1,4,1]$  & $1$ & $D_6$, 8th, $a=-1$\\

$[-2,5,4,-6,10,-12,10 ]$ & $[1, \frac{1}{2},\frac{1}{2},\frac{1}{2}, \frac{1}{2},\frac{1}{2},\frac{1}{2}]$ & $[1,0,1,1,2,1,1]$  & $1$ & $E_6$, 4th, $g=-1$\\

$[9,-1,-11,9,1,-3,-1 ]$ & $[\frac{1}{2},1, \frac{1}{2},\frac{1}{2},\frac{1}{2}, \frac{1}{2},1]$ & $[1,2,3,0,1,0,2]$  & $1$  & $E_6$, 4th*, $g=-3$\\

$[-6,-10,5,8,1,-3,-2]$ & $[\frac{1}{2},\frac{1}{2},\frac{1}{2}, \frac{1}{2},\frac{1}{2},\frac{1}{2},1]$ & $[1,2,2,0,1,1,1 ]$  & $1$ & $E_6$, 5th*, $g=-1$\\

$[-11,5,9,-7,5,5,-7]$ & $\frac{\rho}{2}$ & $\rho$ & $1$ & $E_6$, 7th, $g=-1$\\

$[9,-1,-11,10,-1,-2,-1]$ & $[\frac{1}{2},1,\frac{1}{2},\frac{1}{2},1,1,1]$ & $[1,0,7,0,0,1,0 ]$ & $1$ & $A_1.A_5$, $c=-1$\\

$[-1,-1,-6,4,9,-11,9 ]$ & $[1, 1, \frac{1}{2},\frac{1}{2}, \frac{1}{2},\frac{1}{2}, \frac{1}{2}]$ & $[2,0,2,0,1,3,1 ]$ & $1$ & $D_6$, 8th, $a=-3$\\

$[-2,-12,-1,10,-1,5,-6 ]$ & $[1, \frac{1}{2}, 1, \frac{1}{2}, 1, \frac{1}{2}, \frac{1}{2}]$ & $[1,6,0,1,0,0,4 ]$ & $1$  & $E_6$, 8th, $g=-1$\\

$[-6,-10,5,9,-1,-2,-2 ]$ & $[\frac{1}{2},\frac{1}{2}, \frac{1}{2},\frac{1}{2}, 1,1,1 ]$ & $[4,6,0,0,0,1,1 ]$ & $1$  & $A_6$, 13th*, $b=-1$\\

$[-9,7,7,-9,7,5,-7 ]$ & $\frac{\rho}{2}$ & $\rho$ & $1$ & $E_6$, 9th, $g=-1$\\

$[-2,9,8,-10,8,-10,8 ]$ & $[1,\frac{1}{2},\frac{1}{2},\frac{1}{2},\frac{1}{2},\frac{1}{2},\frac{1}{2} ]$ & $[1,0,1,1,2,1,1]$ & $1$ & $A_4.A_2$, $e=-3$\\

$[-10,-1,8,4,-6,8,-10 ]$ & $[\frac{1}{2},1, \frac{1}{2}, \frac{1}{2},\frac{1}{2}, \frac{1}{2},\frac{1}{2}]$ & $[2,0,1,1,1,2,1 ]$ & $1$ & $E_6$, 6th, $g=-3$\\

$[ -1,-13,-2,11,3,-4,-1]$ & $[1, \frac{1}{2}, 1, \frac{1}{2},\frac{1}{2}, \frac{1}{2},1]$ & $[2,5,0,1,0,0,3 ]$ & $1$ & $E_6$, 8th, $g=-3$\\

$[-10,-1,8,6,-8,6,-8 ]$ & $[\frac{1}{2}, 1, \frac{1}{2},\frac{1}{2}, \frac{1}{2}, \frac{1}{2}, \frac{1}{2}]$ & $[2,0,1,1,1,2,1 ]$ & $1$ & $D_6$, 14th, $g=-1$\\

$[-2,-1,9,-11,9,6,-7 ]$ & $[1,1, \frac{1}{2},\frac{1}{2}, \frac{1}{2}, \frac{1}{2}, \frac{1}{2} ]$ & $[1,3,0,2,0,0,3 ]$ & $1$ & $E_6$, 11th, $g=-1$\\

$[ -12,8,11,-9,7,-2,-2]$ & $[\frac{1}{2},\frac{1}{2}, \frac{1}{2}, \frac{1}{2}, \frac{1}{2},1,1 ]$ & $[1,0,4,0,1,1,1]$  & $1$ & $D_6$, 11th, $a=-3$ \\

$[1,-3,-3,1,11,-13,11 ]$ & $\frac{\rho}{2}$ & $\rho$ & $1$  & $E_6$, 4th, $g=-3$\\

$[-10,10,9,-11,9,-2,-2 ]$ & $[\frac{1}{2},\frac{1}{2}, \frac{1}{2}, \frac{1}{2}, \frac{1}{2},1,1  ]$ & $[1,0,4,0,1,1,1 ]$ & $1$ & $A1.A5$, $c=-3$\\

$[-1,-2,-1,-1,12,-13,11 ]$ & $[1,1,1,1,\frac{1}{2},\frac{1}{2}, \frac{1}{2} ]$ & $[0,1,0,0,0,9,1]$ & $1$ & $D_5.A_1$, $f=-1$\\

$[ -1,-14,-1,13,-2,-1,-1]$ & $[1, \frac{1}{2}, 1, \frac{1}{2}, 1, 1, 1]$ & $[0,10,0,0,1,0,0 ]$ & $1$ & $A_6$, $b=-1$\\

$[-12,-2,11,-1,-1,8,-10 ]$ & $[\frac{1}{2}, 1, \frac{1}{2}, 1, 1, \frac{1}{2}, \frac{1}{2} ]$ & $[8,1,0,0,0,1,6 ]$ & $1$ & $D_6$, 17th, $a=-1$\\

$[-14,-7,12,5,-6,4,-2]$ & $[\frac{1}{2},\frac{1}{2},\frac{1}{2},\frac{1}{2},\frac{1}{2},\frac{1}{2},1 ]$ & $[1,0,2,1,1,1,1 ]$ & $1$& $D_6$, 19th*, $a=-1$\\

$[-1,-2,12,-13,12,-2,-1]$ & $[1,1,\frac{1}{2},\frac{1}{2},\frac{1}{2},1,1 ]$ & $[0,7,0,0,2,0,0 ]$ &$1$ & $A_6$, $b=-3$\\

$[-1,-1,-11,9,-1,9,-11 ]$ & $[1,1,\frac{1}{2},\frac{1}{2},1,\frac{1}{2},\frac{1}{2} ]$ & $[6,2,0,0,0,1,5 ]$ & $1$&$D_6$, 17th, $a=-3$\\

$[-1,-2,-1,9,-11,13,-14 ]$ & $[1,1,1,\frac{1}{2},\frac{1}{2},\frac{1}{2},\frac{1}{2} ]$ & $[ 0,1,0,1,0,6,1]$ & $1$&$E_6$, 13th, $g=-3$\\

$[-13,-2,12,-1,6,-8,-1]$ & $[\frac{1}{2},1,\frac{1}{2},1,\frac{1}{2},\frac{1}{2},1 ]$ & $[7,1,0,0,0,2,4 ]$ &$1$  & $D_6$, 20th, $a=-1$\\

$[-1,-2,-1,11,-13,11,-12]$ & $[1,1,1,\frac{1}{2},\frac{1}{2},\frac{1}{2},\frac{1}{2} ]$ & $[0,1,0,1,0,6,1 ]$ &  $1$ & $D_5.A_1$, $f=-3$\\

$[-1,-1,-12,10,7,-9,-1 ]$ & $[1,1,\frac{1}{2},\frac{1}{2},\frac{1}{2},\frac{1}{2},1 ]$ & $[ 5,2,0,0,0,2,3]$ &  $1$& $D_6$, 20th, $a=-3$\\

$[-14,-2,13,4,-6,-1,-1]$ & $[\frac{1}{2},1,\frac{1}{2},\frac{1}{2},\frac{1}{2},1,1 ]$ & $[6,1,0,0,0,3,2 ]$ &  $1$ & $D_6$, 21st, $a=-1$\\

$[-2,8,-1,-9,7,10,-12 ]$ & $[1,\frac{1}{2},1,\frac{1}{2},\frac{1}{2},\frac{1}{2},\frac{1}{2} ]$ & $[4,3,0,0,0,1,4 ]$ & $1$ &$D_6$, 17th, $a=-5$\\

$[-1,7,7,-9,-1,13,-15 ]$ & $[1,\frac{1}{2},\frac{1}{2},\frac{1}{2},1,\frac{1}{2},\frac{1}{2} ]$ & $[0,2,0,1,0,4,2 ]$ & $1$ &$E_6$, 13th, $g=-5$\\

$[-17,1,15,-3,1,-3,1 ]$ & $\frac{\rho}{2}$ & $\rho$ & $1$ & $D_6$, 22nd, $a=-1$\\

$[11,-2,-12,11,-13,11,-2 ]$ & $[\frac{1}{2},1,\frac{1}{2},\frac{1}{2},\frac{1}{2},\frac{1}{2},1 ]$ & $[0,4,1,0,2,0,1 ]$ & $1$ & $A_6$, $b=-5$\\

$[-1,8,7,-9,11,-13,-2 ]$ & $[1,\frac{1}{2},\frac{1}{2},\frac{1}{2},\frac{1}{2},\frac{1}{2},1 ]$ & $[0,1,0,2,0,3,1]$  & $1$ &$E_6$, 15th*, $g=-5$\\

$[-1,10,9,-11,9,-11,-2 ]$ & $[1,\frac{1}{2},\frac{1}{2},\frac{1}{2},\frac{1}{2},\frac{1}{2},1 ]$ & $[0,1,0,2,0,3,1]$ & $1$ & $D_5.A_1$, $f=-5$
\end{tabular}
\label{table-E7-scattered-part-1}
\end{table}

\begin{table}
\centering
\caption{The scattered representations in $\widehat{E}_7^{\mathrm{d}}$ (part two)}
\begin{tabular}{r|c|c|c|c}
$s\rho$ &   $\lambda$   & spin LKT & mult  & string limits\\
\hline

$[-17,-1,15,-1,-1,-1,-1]$ & $[\frac{1}{2},1,\frac{1}{2},1,1,1,1]$ & $[13,0,1,0,0,0,0]$ &$1$ & $D_6$,  $a=-1$\\

$[4,-6,-6,4,-1,14,-16]$ & $[\frac{1}{2},\frac{1}{2},\frac{1}{2},\frac{1}{2},1,\frac{1}{2},\frac{1}{2}]$ & $[ 1,2,0,1,0,3,2]$ & $1$ & $E_6$, 13th, $g=-7$\\

$[-1,-1,-16,14,-1,-1,-1]$ & $[1, 1, \frac{1}{2}, \frac{1}{2}, 1, 1, 1]$ & $[10,0,2,0,0,0,0 ]$ &$1$ & $D_6$, $a=-3$\\

$[4,-7,-6,5,12,-14,-2 ]$ & $[\frac{1}{2}, \frac{1}{2}, \frac{1}{2}, \frac{1}{2}, \frac{1}{2}, \frac{1}{2}, 1 ]$ & $[1,1,0,2,0,2,1]$ & $1$&$E_6$, 15th*, $g=-7$\\

$[-11,9,9,-11,9,-11,9 ]$ & $\frac{\rho}{2}$ & $\rho$ & $1$ & $A_6$, $b=-7$\\

$[-3,-1,1,-3,1,15,-17 ]$ & $[\frac{1}{2}, 1, \frac{1}{2}, \frac{1}{2}, \frac{1}{2}, \frac{1}{2}, \frac{1}{2} ]$ & $[2,2,0,1,0,2,2]$ & $1$ & $E_6$, 13th, $g=-9$\\

$[-1,13,-1,-15,13,-1,-1 ]$ & $[1,\frac{1}{2},1,\frac{1}{2},\frac{1}{2},1,1]$ & $[7,0,3,0,0,0,0 ]$ &$1$ & $D_6$,  $a=-5$\\

$[-2,-1,-1,-1,-1,16,-17 ]$ & $[ 1,1,1,1,1, \frac{1}{2},\frac{1}{2}]$ & $[1,0,0,0,0,0,15 ]$ & $1$ & $E_6$, $g=-1$\\

$[-1,-1,-2,-1,15,-16,-1 ]$ & $[ 1, 1, 1, 1, \frac{1}{2},\frac{1}{2}, 1]$ & $[2,0,0,0,0,0,13 ]$ &$1$ & $E_6$, $g=-3$ \\

$[-1,-13,-1,11,-13,11,-1 ]$ & $[1, \frac{1}{2}, 1, \frac{1}{2}, \frac{1}{2}, \frac{1}{2}, 1 ]$ & $[5,0,3,0,0,1,0]$ & $1$ & $D_6$,  $a=-7$\\

$[-1,-1,-1,13,-15,-1,-1]$ & $[1, 1, 1, \frac{1}{2},\frac{1}{2}, 1, 1]$ & $[ 3,0,0,0,0,0,11]$ &$1$ & $E_6$, $g=-5$\\

$[-1,-1,9,-11,9,-11,9 ]$ & $[1, 1, \frac{1}{2}, \frac{1}{2}, \frac{1}{2}, \frac{1}{2}, \frac{1}{2}]$ & $[3,0,3,0,0,2,0]$ & $1$&$D_6$,  $a=-9$\\

$[-1,11,12,-13,-1,-2,-1 ]$ & $[1, \frac{1}{2}, \frac{1}{2}, \frac{1}{2}, 1, 1, 1 ]$ & $[4,0,0,0,0,0,9]$ & $1$& $E_6$, $g=-7$\\

$[10,-10,-11,9,-1,-2,-2]$ & $[\frac{1}{2}, \frac{1}{2}, \frac{1}{2}, \frac{1}{2}, 1, 1, 1 ]$ & $[4,0,0,0,0,1,7]$ & $1$ & $E_6$,  $g=-9$\\

$[6,-1,-8,6,-8,6,-8 ]$ & $[\frac{1}{2}, 1, \frac{1}{2}, \frac{1}{2}, \frac{1}{2}, \frac{1}{2}, \frac{1}{2}]$ & $[2,0,2,1,0,2,0]$ & $1$ &$D_6$, $a=-11$\\

$[-9,-1,7,-8,6,-1,-2 ]$ & $[\frac{1}{2}, 1, \frac{1}{2}, \frac{1}{2}, \frac{1}{2}, 1, 1]$ & $[4,0,0,0,0,2,5]$ & $1$ & $E_6$, $g=-11$\\

$[-1,-2,-6,4,-5,4,-2 ]$ & $[1,1, \frac{1}{2}, \frac{1}{2}, \frac{1}{2}, \frac{1}{2}, 1]$ & $[4,0,0,0,0,3,3]$ & $1$ & $E_6$, $g=-13$\\

$[-5,3,3,-5,3,-5,-1 ]$ & $[\frac{1}{2}, \frac{1}{2}, \frac{1}{2}, \frac{1}{2}, \frac{1}{2}, \frac{1}{2}, 1]$ & $[1,0,1,2,0,2,0]$ & $1$ &$D_6$, $a=-13$\\

$[-1,1,-1,-3,1,-3,1 ]$ & $\frac{\rho}{2}$ & $\rho$ & $1$ & \\

$[-1,1,-1,-3,1,-3,1 ]$ & $[1, \frac{1}{2}, 1, \frac{1}{2}, \frac{1}{2}, \frac{1}{2}, \frac{1}{2} ]$ & $[4,0,0,0,0,4,1]$ & $1$& $E_6$, $g=-15$\\

$-\rho$ & $[\frac{1}{2}, \frac{1}{2}, \frac{1}{2}, \frac{1}{2}, \frac{1}{2}, \frac{1}{2}, 1]$ & $[1,0,1,2,0,2,0]$ & $1$ &\\

$-\rho$ & $\rho$ & [0,0,0,0,0,0,0] & $1$ &
\end{tabular}
\label{table-E7-scattered-part-2}
\end{table}

Note that there are three scattered  representations of complex $E_7$ which can not be realized as string limits. Namely, the 19th, 21st and 22nd entries of Table \ref{table-E7-scattered-part-2}. They are all unipotent representations. We will look at them more carefully in the next section.

By the way, we note that for complex $E_6$, only the two \emph{spherical} scattered representations can \emph{not} be realized as string limits: the model representation \cite{Mc}, and the trivial representation.

\section{Unipotent representations} \label{sec-unip}


Let $G$ be a reductive Lie group over a local field $F$. In \cite{A}, Arthur studied the homomorphisms
$$\Phi : W_F \times SL(2, \mathbb{C}) \to\ ^LG,$$
where $W_F$ is the Weil group for a local field $F$ and $^LG$ is the Langlands dual group of $G$. For each $\Phi$, he
conjectured that there should be a packet of representations related to certain
automorphic representations of locally symmetric spaces.
We now focus on the case when $F = \mathbb{C}$. Suppose $\Phi|_{W_F} = Triv$, and the image of $d\Phi$
is a special nilpotent orbit (in the sense of Lusztig) in $^L\mathfrak{g}$, then the packet of representations
corresponding to $\Phi$ are called {\it special unipotent}.


Special unipotent representations are studied in \cite{BV} in full detail. The main result is that they satisfy the properties conjectured by Arthur.
For classical groups, the classification of the unitary dual given in \cite{V1} and \cite{B1} verified that all such representations are unitary.
It is also easy to see that the same holds for $G_2$ (see \cite{Du} for a full classification of the unitary dual of $G_2$).
However, their unitarity are not completely known for the other exceptional groups.

One may conjecture that the special unipotent representations form the building blocks of the unitary spectrum of $G$, that is, all elements in $\widehat{G}$ can be obtained from them
by unitary induction, along with some continuous deformations on the induced modules (the complementary series).
However, this is obviously not true for groups outside Type $A$ --
for example, the oscillator representations in $Sp(2n,\mathbb{C})$ are not unitarily induced from any representation on a proper
Levi component, and its infinitesimal character is not of the form $^\vee h/2$ as required by \cite{BV}.
For this reason, \cite{B1} and \cite{B2} considered a larger class of representations called {\it unipotent representations}
for all classical groups (in fact, the definition given in \cite[Chapter 11]{B1} can be extended to all simple Lie groups).

Following \cite{B2}, for each
nilpotent orbit $\mathcal{O} \subset \mathfrak{g}$, a {\bf unipotent representation} $\Pi$ attached to $\mathcal{O}$ is an irreducible representation with a fixed infintesimal character $(\lambda_{\mathcal{O}},\lambda_{\mathcal{O}})$ possessing the following properties:
\begin{enumerate}
\item The associated variety of $\Pi$ is equal to the closure of ${\mathcal{O}}$;
\item $Ann(\Pi) \subset U(\mathfrak{g})$ is the maximal primitive ideal
  $\mathcal{I}_{\lambda_{\mathcal{O}}}$ with infinitesimal character $(\lambda_{\mathcal{O}},\lambda_{\mathcal{O}})$;
\item {$|\{\Pi : Ann(\Pi)=\mathcal{I}_{\lambda_{\mathcal{O}}}\}| = | \widehat{A(\mathcal{O})}|,$ where $A({\mathcal{O}})$ is the component group of the centralizer of an $e\in \mathcal{O}$};
\item $\Pi$ is unitary.
\end{enumerate}
From now on, we say $\Pi$ is unipotent if all the above conditions hold. In the case when a special orbit $\mathcal{O}$ is {\it stably trivial}, i.e. $A(\mathcal{O})$ is isomorphic to the {\it Lusztig's quotient} $\overline{A}(\mathcal{O})$, the above definition shall match with that of special unipotent representation in \cite{BV}.

\medskip
We apply the above formulation of unipotent representation for $E_7$. In \cite{Mc}, it was shown that the representations listed below
satisfy (1) and (2) for the specified $\mathcal{O}$ (using the Bala-Carter notation  \cite[Chapter 8]{CM}) and $\lambda_{\mathcal{O}}$:

\smallskip
\noindent \underline{$\mathcal{O} = 4A_1, \lambda_{\mathcal{O}} = \rho/2:$}
\begin{itemize}
\item[(i)]  $J(\frac{\rho}{2}, \frac{\rho}{2})$ (the 1st entry of \cite[Table 10]{Mc});
\item[(ii)]  $J(\frac{\rho}{2}, [\frac{1}{2}, -\frac{1}{2}, \frac{1}{2}, \frac{3}{2}, -\frac{1}{2}, \frac{3}{2}, -\frac{1}{2}])$ (cf. $Q'$ in the proof of \cite[Theorem 2.2]{Mc});
\end{itemize}

\smallskip
\noindent \underline{$\mathcal{O} = (3A_1)', \lambda_{\mathcal{O}} = [\frac{1}{2}, \frac{1}{2},\frac{1}{2}, \frac{1}{2}, \frac{1}{2}, \frac{1}{2}, 1]:$}
\begin{itemize}
\item[(iii)]  $J(\lambda_{\mathcal{O}}, \lambda_{\mathcal{O}})$ (the 2nd entry of \cite[Table 10]{Mc});
\end{itemize}

\smallskip
\noindent \underline{$\mathcal{O} = (3A_1)'', \lambda_{\mathcal{O}} = [1, \frac{1}{2}, 1, \frac{1}{2}, \frac{1}{2}, \frac{1}{2}, \frac{1}{2}]:$}
\begin{itemize}
\item[(iv)]  $J(\lambda_{\mathcal{O}}, \lambda_{\mathcal{O}})$ (the 3rd entry of \cite[Table 10]{Mc});
\end{itemize}

\smallskip
\noindent \underline{$\mathcal{O} = 2A_1, \lambda_{\mathcal{O}} = [1,1,1,0,1,0, 1]:$}
\begin{itemize}
\item[(v)]   $J(\lambda_{\mathcal{O}}, \lambda_{\mathcal{O}})$ (the 4th entry of \cite[Table 10]{Mc});
\end{itemize}

\smallskip
\noindent \underline{$\mathcal{O} = A_1, \lambda_{\mathcal{O}} = [1,1,1,0,1,1, 1]:$}
\begin{itemize}
\item[(vi)]   $J(\lambda_{\mathcal{O}}, \lambda_{\mathcal{O}})$ (the  minimal representation);
\end{itemize}

\smallskip
\noindent \underline{$\mathcal{O} = 0, \lambda_{\mathcal{O}} = \rho:$}
\begin{itemize}
\item[(vii)]  $J(\rho, \rho)$, the trivial representation.
\end{itemize}
Moreover, all of these representations are unitary (David Vogan helped us with the representation (v), which used about 290 GB memory in \texttt{atlas}). So (4) also holds for these representations, and we are left to check whether (3) holds or not.

As suggested by Jeffery Adams, by looking the table in
\begin{center}
\texttt{liegroups.org/tables/unipotentExceptional/E7/E7\_ad/E7\_ad\_summary.txt}
\end{center}
on the centralizer of any nilpotent element $e \in \mathcal{O}$, one expects
the number of representations satisfying (3) to be:
$$4A_1 : 2; \quad \quad (3A_1)' : 1; \quad  \quad (3A_1)'' : 2; \quad  \quad 2A_1 : 1;\quad \quad A_1 : 1;\quad \quad 0 : 1.$$
Note that the for the {\it model orbit} $4A_1$, the fundamental group should be $\mathbb{Z}/2\mathbb{Z}$ instead of $1$ given by the table in \cite{CM}. The same observation is also made in Section 4.6.1 and Remark 4.13 of \cite{C}.

These numbers match with our list (i) -- (vii) above, except for $((3A_1)'', [1, \frac{1}{2}, 1, \frac{1}{2}, \frac{1}{2}, \frac{1}{2}, \frac{1}{2}])$ where we are off by one representation. Indeed, (iv) and $\Pi_{20}$, the 20th entry of Table \ref{table-E7-scattered-part-2}, both satisfy (1) -- (4) above. More precisely, (iv) $= \mathrm{Ind}_{E_6 \times \mathbb{C}^*}^{E_7}(\mathrm{triv} \boxtimes |\det|^{1/2})$ (mid-point of complementary series) and $\Pi_{20} = \mathrm{Ind}_{E_6 \times \mathbb{C}^*}^{E_7}(\mathrm{triv} \boxtimes \det)$. As a consequence, one can conclude that all the representations (i) -- (vii) along with $\Pi_{20}$ are unipotent.

We look more closely to the case of $\mathcal{O} = (3A_1)''$. It is a speical Richardson orbit, which is induced from the trivial orbit of $E_6$.
The 2-fold cover $\widetilde{\mathcal{O}}$ of $\mathcal{O}$ is {\it birationally rigid}, i.e. it cannot be birationally induced from any proper Levi subgroup. We make the following remarks:
\begin{itemize}
\item $\mathcal{O}$ is {\bf not} stably trivial, and it is the only non-stably trivial special orbit in our list above. More precisely, $A(\mathcal{O}) = \mathbb{Z}/2\mathbb{Z}$ but $\overline{A}(\mathcal{O}) = \{e\}$. So the representations (iv) and $\Pi_{20}$ are unipotent but not special unipotent. Indeed, the parameter of the special unipotent representation
is $h(^\vee \mathcal{O})/2 = [1,0,1,1,0,1,0]$ (instead of $[1, \frac{1}{2}, 1, \frac{1}{2}, \frac{1}{2}, \frac{1}{2}, \frac{1}{2}]$), and is equal to $\mathrm{Ind}_{E_6 \times \mathbb{C}^*}^{E_7}(\mathrm{triv} \boxtimes \mathrm{triv})$.
\item In \cite{LMM} and \cite{MM}, another definition of unipotent representation is given in terms of the {\it orbit philosophy}. More explicitly, they assign a unipotent parameter corresponding to each $G$-equivariant cover of a nilpotent orbit. Since $\mathcal{O}$ is Richardson from a Levi of Type $E_6$, the parameter corresponding to $\widetilde{\mathcal{O}}$ is given by $\rho(E_6)$, half the sum of roots in the subroot system of $E_6$, which in turn can be conjugated to the special unipotent parameter $[1,0,1,1,0,1,0]$. On the other hand, since $\widetilde{\mathcal{O}}$ is birationally rigid, one needs to `shift' $\rho(E_6)$ to obtain the parameter corresponding to $\widetilde{\mathcal{O}}$. By the calculations in \cite[Section 4.3.2]{MM},
it is equal to the unipotent parameter $[1, \frac{1}{2}, 1, \frac{1}{2}, \frac{1}{2}, \frac{1}{2}, \frac{1}{2}]$ in our list.
\end{itemize}

To finish our study on unipotent representations, we check whether (i) -- (vii) and $\Pi_{20}$
are in $\widehat{G}^d$ or not. Indeed, (i) and (iv) have zero Dirac cohomology, since neither $\rho$ nor $[2, 1, 2, 1, 1, 1, 1]$ lies in the root lattice of $E_7$ (recall Equation \eqref{Lambda-s}),
while (v) and (vi) also have zero Dirac cohomology since $\lambda$ is singular. The remaining four representations are all Dirac series. Indeed, (ii), $\Pi_{20}$, (iii), (vii) are the 19th -- 22nd entries of Table \ref{table-E7-scattered-part-2}, respectively.

These discussions lead us to Corollary \ref{cor-main}.

\section{Calculations of Spin Lowest $K$-types}
We mention some calculations on the spin lowest $K$-types $\sigma$ in Table \ref{table-E7-scattered-part-2}. It is easy to check that all of them satisfy $\|\sigma\|_{\mathrm{spin}} = \|\Lambda\| = \|2\lambda\|$.

Firstly, let us handle the cases that $J(\lambda,-s\lambda) \in \widehat{G}^d$ can be realized as a string limit. In such a case, it suffices to
check that
\begin{equation}\label{string-limit-spinLKT}
[V(\sigma) : J(\lambda,-s\lambda)] > 0.
\end{equation}
Indeed, the branching result in \cite{BDW} implies that $\sigma$ must be the unique spin lowest and it must occur exactly once. Usually, one can use the $\texttt{atlas}$ command $\texttt{branch\_irr}$ to check \eqref{string-limit-spinLKT} directly. However, for a few representations, such a direct calculation in \texttt{atlas} would be time-consuming. Let us mention some ways that are more effective.

\subsection*{$\bullet$ $E_6$, $g = -15$ (the 20th entry of Table \ref{table-E7-scattered-part-2})}

As discussed in the previous section, this is the unipotent representation $\Pi_{20} = \mathrm{Ind}_{E_6 \times \mathbb{C}^*}^{E_7}(\mathrm{triv} \boxtimes \det)$ having LKT $[0,0,0,0,0,0,1]$ of height $27$, while the $K$-type $[4, 0, 0, 0, 0, 4, 1]$ has height $371$. Instead of doing \texttt{branch\_irr(p, 371)}, we do \texttt{branch\_irr(p, 235)} and find that the $K$-type $[0, 0, 0, 0, 0, 4, 1]$ shows up in \texttt{p} exactly once. Using Vogan pencil, we conclude that the $K$-type
$$
\sigma:=[0, 0, 0, 0, 0, 4, 1]+ 4 \beta=[4, 0, 0, 0, 0, 4, 1]
$$
shows up in \texttt{p}, as desired.

It is interesting to note that if we look at the representations $E_6$, $g=-1, -3, \dots, -15$ consecutively, the spin LKTs still demonstrate a vague pattern, whilst the pattern is described explicitly in \cite{DW2} for all classical groups.

There are cases that we can not move \emph{down} from the desired spin LKT along the Vogan pencil.
Then the following clever way (suggested by David Vogan) could be helpful to get the multiplicity of $\sigma$ in $ J(\lambda,-s\lambda)$.

\subsection*{$\bullet$ $D_6$, $a = -13$ (the 18th entry of Table \ref{table-E7-scattered-part-2})} Here $J(\lambda,-s\lambda) = {\rm Ind}_{D_6 \times \mathbb{C}^*}^{E_7}(\mathrm{triv} \boxtimes \det^2)$
with lowest $K$-type  $2\beta=[2,0,0,0,0,0,0]$. As mentioned above,
it suffices to check
$$
[V([1,0,1,2,0,2,0]): {\rm Ind}_{D_6 \times \mathbb{C}^*}^{E_7}(\mathrm{triv} \boxtimes \det^2)] > 0.
$$
Using Frobenius reciprocity, it is equivalent to show that
\begin{equation} \label{eq-branch1}
[\mathrm{triv} \boxtimes \mathbb{C}_{2\beta} : V([1,0,1,2,0,2,0])|_{D_6 \times \mathbb{C}^*}] >0.
\end{equation}
The branching $V([1,0,1,2,0,2,0])|_{D_6 \times \mathbb{C}^*}$ can be done as follows:
$V([1,0,1,2,0,2,0])$ defines a finite-dimensional $(\mathfrak{e}_7, D_6 \times A_1)$-module for the quarterionic real form \texttt{E7\_q},
so the $K$-types of this finite-dimensional module are precisely the branching of $V([1,0,1,2,0,2,0])$ to $D_6 \times A_1$.

The following \texttt{atlas} command gives the finite-dimensional module $V([1,0,1,2,0,2,0])$ and its dimension.
\begin{verbatim}
set G=E7_q
set p=parameter(x_open(G),[1,0,1,2,0,2,0]+rho(G),[1,0,1,2,0,2,0]+rho(G))
dimension(p)
Value: 2399133156669849600
\end{verbatim}
Branching of the module \texttt{p} into $D_6 \times A_1$-modules could cost a lot of time. However, it is easy to see that it starts as follows:
\begin{verbatim}
print_branch_irr_long(p,KGB(G,62),40)
m  x    lambda                            hw                       dim  height
1  20   [ 0, 0, 0, 0, 0, 0, 1 ]/1         [ 1, 0, 0, 0, 0, 0, 0 ]  3    27
1  337  [  1,  0, -1,  1, -1,  1,  0 ]/1  [ 2, 0, 0, 0, 0, 0, 0 ]  5    38
\end{verbatim}
with the last 6 coordinates of \texttt{hw} are the highest weights of the $D_6$-module in the decomposition.
Therefore,
$$
V([1,0,1,2,0,2,0])|_{D_6 \times A_1} = (\mathrm{triv} \boxtimes \mathbb{C}^3) \oplus (\mathrm{triv} \boxtimes \mathbb{C}^5) \oplus \cdots.
$$
Here $\bbC^3$ means the 3-dimensional irreducible representation of $A_1$.
Thus the restriction
$$
V([1,0,1,2,0,2,0])|_{D_6 \times \mathbb{C}^*}
$$
decomposes as
$$
\left(\mathrm{triv} \boxtimes (\mathbb{C}_{\beta} \oplus \mathbb{C}_0 \oplus \mathbb{C}_{-\beta})\right)  \oplus \left(\mathrm{triv} \boxtimes  (\mathbb{C}_{2\beta} \oplus \mathbb{C}_{\beta} \oplus \mathbb{C}_0 \oplus \mathbb{C}_{-\beta} \oplus \mathbb{C}_{-2\beta})\right) \oplus \cdots
$$
and \eqref{eq-branch1} has multiplicity $ > 0$, as desired.

We have also handled the  12th and 15th entry of Table \ref{table-E7-scattered-part-2} in this way.

Secondly, let us determine the spin lowest $K$-types of the two unipotent representations.

\subsection*{$\bullet$ $J(\frac{\rho}{2}, [\frac{1}{2}, -\frac{1}{2}, \frac{1}{2}, \frac{3}{2}, -\frac{1}{2}, \frac{3}{2}, -\frac{1}{2}])$ (the 19th entry of Table \ref{table-E7-scattered-part-2})}
As mentioned in the previous section, this is the unipotent representation $Q'$ appearing in the proof of \cite[Theorem 2.2]{Mc},
whose $K$-types have highest weights in the weight lattice \emph{but not in the root lattice}, all appearing with multiplicity one.
In particular $\sigma = \rho$ appears in $Q'$, and this is the only $\sigma$ that satisfies $\{\sigma - \rho\} + \rho = 2\lambda = \rho $.

\subsection*{$\bullet$$J(\lambda, \lambda)$, $\lambda=[\frac{1}{2}, \frac{1}{2},\frac{1}{2}, \frac{1}{2}, \frac{1}{2}, \frac{1}{2}, 1]$ (the 21st entry of Table \ref{table-E7-scattered-part-2})}
According to \cite[Table 10]{Mc}, it is conjectured that this unipotent representation would have $K$-types equal to $[a,0,b,c,0,d,0]$ with multiplicity one, where $a, b, c, d$ are non-negative integers.

This can be verified by \cite{MM}, which states that this representation is \emph{the} quantization model for the birationally rigid orbit
$\mathcal{O} = 3A_1'$, and hence its $K$-spectrum is equal to the ring of regular functions of $\mathcal{O}$. Then the conjecture follows from the
results of \cite{C}. It is then easy to check that $[1,0,1,2,0,2,0]$ is the \emph{unique} spin lowest $K$-type in $J(\lambda, \lambda)$.

It is also worth noting that this unipotent module, along with the module with string limit $D_6$, $a = -13$ (the 18th entry of Table \ref{table-E7-scattered-part-2}), are the only composition
factors of the end point of the complementary series in \cite[Equation (5.29)]{BP}, each occurring once. Hence both modules are unitary.

\medskip

\centerline{\scshape Acknowledgments}

 Our sincere gratitude is expressed to Professor Vogan for testing the unitarity of the representation (v) in Section \ref{sec-unip}, and to Professor Adams for sharing his knowledge on the model orbit $4A_1$ in Section \ref{sec-unip} with us. We also thank Dr. Luan for helping us considering the branching \eqref{eq-branch1}. We sincerely thank an anonymous referee for his/her great patience in reading this manuscript and offering us excellent suggestions.

\medskip

\centerline{\scshape Funding}
Dong is supported by the National Natural Science Foundation of China (grant 12171344). Wong is supported by the National Natural Science Foundation of China (grant 11901491) and the Presidential Fund of CUHK(SZ).

\end{document}